\documentclass[a4paper,12pt,reqno]{amsart}
\IfFileExists{lmodern.sty}{}{\pdfmapfile{+cm.map}\pdfmapfile{+symbols.map}\pdfmapfile{+cmextra.map}}

\IfFileExists{lmodern.sty}{\usepackage[T1]{fontenc}\usepackage{lmodern}}{}
\usepackage[utf8]{inputenc}
\usepackage[left=2.2cm,right=2.2cm,top=2.2cm,bottom=2.7cm]{geometry}
\usepackage{amsthm,amsmath,amsfonts,amssymb,mathtools}
\usepackage{graphicx}
\usepackage{placeins}
\usepackage[numbers]{natbib}
\usepackage[colorlinks,hyperfootnotes=true,hyperindex,linkcolor=blue,citecolor=blue,urlcolor=blue]{hyperref}
\usepackage{xcolor}

\allowdisplaybreaks
\numberwithin{equation}{section}

\theoremstyle{plain}
\newtheorem{theorem}{Theorem}[section]
\newtheorem{lemma}[theorem]{Lemma}
\newtheorem{proposition}[theorem]{Proposition}
\newtheorem{corollary}[theorem]{Corollary}

\theoremstyle{definition}

\theoremstyle{remark}
\newtheorem{remark}[theorem]{Remark}

\newcommand{\R}{\mathbb{R}}
\renewcommand{\P}{\mathbb{P}}
\newcommand{\E}{\mathbb{E}}
\newcommand{\Var}{\operatorname{Var}}

\newcommand{\abs}[1]{\left|#1\right|}

\newcommand{\1}{\mathbf{1}}

\title[Stein's method for the annealed dilute Curie--Weiss model]
{Stein's method for the annealed dilute Curie--Weiss model:\ Berry--Esseen bounds and a critical crossover}
\author{Peter Eichelsbacher}
\address{Faculty of Mathematics, Ruhr University Bochum, Universit\"atsstra\ss e 150, 44780 Bochum, Germany}
\email{peter.eichelsbacher@rub.de}
\date{November 2026}

\begin{document}
\begin{abstract}
We investigate fluctuations of the magnetisation in the annealed dilute Curie--Weiss model by Stein's method. An exact representation identifies the annealed spin law with a classical Curie--Weiss model at an effective, size-dependent inverse temperature. A uniform estimate for this parameter permits an $O(N^{-1/2})$ Berry--Esseen bound in Kolmogorov distance, together with concentration inequalities and a Cram\'er-type normal-tail approximation, for fixed inverse temperature $0<\beta<1$ and arbitrary fixed external field under the full dense-regime assumption $p_NN\to\infty$. At nominal criticality $\beta=1$ and zero field, we identify three dilution regimes separated by the threshold $p_N\asymp N^{-3/4}$. In the crossover regime $p_N^4N^3\to c^{-2}$ with $c>0$, the rescaled magnetisation has a non-Gaussian limit with density proportional to $\exp\{cx^2/48-x^4/12\}$. The corresponding Stein drift is nonmonotone, so the usual monotonicity-based estimates for indicator test functions do not apply. We derive direct Stein-factor bounds for this double-well family and prove a quantitative Kolmogorov estimate of order $N^{-1/2}$ plus the discrepancy of the effective critical parameter from its limit. We also determine the optimal uniform derivative factor for indicator Stein solutions and identify its transition as the double-well parameter increases.
\end{abstract}
\maketitle

\section{Introduction}\label{sec:introduction}

We consider $N$ spins $\sigma=(\sigma_1,\ldots,\sigma_N)$ in
$\Omega_N=\{-1,1\}^N$. First choose a random directed graph on
$[N]$: independently for each ordered pair $(i,j)\in[N]^2$,
an edge from $i$ to $j$ is present with probability $p_N\in(0,1]$.
Let $\varepsilon_{ij}\in\{0,1\}$ indicate its presence; self-loops
are allowed. Only pairs with $\varepsilon_{ij}=1$ contribute to
the interaction sum below. The random graph is the \emph{disorder}
in the model. We use the graph convention of \cite{ADS26}. Write
\[
 M_N(\sigma)=\sum_{i=1}^N\sigma_i
\]
and consider the dilute Curie--Weiss Hamiltonian
\[
 H_{N,\varepsilon,\beta,h}(\sigma)
 =-\frac{\beta}{2p_NN}\sum_{i,j=1}^N\varepsilon_{ij}\sigma_i\sigma_j-hM_N(\sigma).
\]
Here $\beta>0$ is the \emph{inverse temperature}, and $h\in\R$ is the \emph{external magnetic field}. In the high-temperature part of the paper we assume $0<\beta<1$; we also study the critical value $\beta=1$, at $h=0$.
For a fixed graph $\varepsilon$, the Gibbs measure assigns
weight proportional to $e^{-H_{N,\varepsilon,\beta,h}(\sigma)}$;
this is the \emph{quenched} model. We instead average this
weight over the graph before normalising. The resulting
\emph{annealed} measure is
\[
 \mu_{N,p_N,\beta,h}(\sigma)
 =\frac{\E_\varepsilon e^{-H_{N,\varepsilon,\beta,h}(\sigma)}}
 {\sum_{\tau\in\Omega_N}\E_\varepsilon
     e^{-H_{N,\varepsilon,\beta,h}(\tau)}}.
\]
Here $\E_\varepsilon$ averages only the graph variables;
probabilities and expectations below, unless specified otherwise,
refer to the annealed spin measure. Thus the random graph is
integrated out, leaving a deterministic probability law on
$\Omega_N$.

The annealed law has a useful exact form. Set
$t_N=\beta/(2p_NN)$. For $s\in\{-1,1\}$, independence of the
edge indicators gives the single-edge identity
$\E_\varepsilon e^{t_N\varepsilon_{ij}s}=1-p_N+p_Ne^{t_Ns}$.
Accordingly, put
\[
 A_N=1-p_N+p_Ne^{t_N},\qquad
 B_N=1-p_N+p_Ne^{-t_N}.
\]
Among the $N^2$ ordered pairs, let $N_+$ and $N_-$ denote
the numbers for which $\sigma_i\sigma_j=1$ and $-1$, respectively.
As $N_++N_-=N^2$ and
$N_+-N_-=\sum_{i,j}\sigma_i\sigma_j=M_N(\sigma)^2$,
\[
 N_+=\frac{N^2+M_N(\sigma)^2}{2},\qquad
 N_-=\frac{N^2-M_N(\sigma)^2}{2}.
\]
Independence of the edges therefore gives
\[
 \E_\varepsilon e^{-H_{N,\varepsilon,\beta,h}(\sigma)}
 =A_N^{N_+}B_N^{N_-}e^{hM_N(\sigma)}
 =(A_NB_N)^{N^2/2}
 \left(\frac{A_N}{B_N}\right)^{M_N(\sigma)^2/2}
 e^{hM_N(\sigma)}.
\]
This is the same exact reduction as in \cite[(8)--(11)]{ADS26},
written here by counting the two types of spin pairs.
The first factor does not depend on $\sigma$ and cancels
under normalisation. Hence the annealed model is exactly
a classical Curie--Weiss model,
\begin{equation}\label{eq:intro-exact-CW}
 \mu_{N,p_N,\beta,h}(\sigma)
 =\frac1{\widetilde Z_N}
 \exp\left\{\frac{\beta_N}{2N}M_N(\sigma)^2+hM_N(\sigma)\right\},
\end{equation}
where the parameter $\beta_N$, called the \emph{effective inverse
temperature}, is given explicitly by
\begin{equation}\label{eq:intro-betaN}
 \beta_N=N\log\frac{1-p_N+p_Ne^{\beta/(2p_NN)}}
 {1-p_N+p_Ne^{-\beta/(2p_NN)}}.
\end{equation}
There is no approximation in this reduction; all effects of
annealed dilution enter through $\beta_N$.
For general $b>0$ and $u\in\R$, we write $CW(b,u)$ for the
classical Curie--Weiss law on $\Omega_N$ given by the right-hand side of
eq.~\eqref{eq:intro-exact-CW} with $\beta_N$ and $h$ replaced by
$b$ and $u$, respectively. Its normalising constant is then understood
to depend on $N,b,u$.

Our basic dilution assumption is the \emph{dense regime}
$p_NN\to\infty$: the expected number $p_NN$ of outgoing
edges per vertex diverges. This is the condition assumed in
our high-temperature results, not a claim that smaller degrees
are mathematically impossible.
Apostel, D\"oring and Schubert \cite{ADS26} established
a broad collection of high-temperature results---including
Berry--Esseen bounds, Cram\'er corrections, concentration,
moderate deviations and mod-Gaussian convergence---under
\[
 0<\beta<1,\qquad h\in\R,\qquad p_N^3N^2\to\infty.
\]
Their cumulant method requires the stronger condition
$p_N^3N^2\to\infty$. Our Stein approach gives Berry--Esseen,
concentration and classical Cram\'er-type estimates under
$p_NN\to\infty$, and also treats $\beta=1$, $h=0$.

The classical Curie--Weiss fluctuation theory goes back to
Ellis and Newman \cite{EN78a,EN78b}. Bovier and Gayrard
\cite{BG93} studied random Bernoulli couplings in the dense
regime. For Stein approximation in mean-field models see
\cite{EL10,CS11}; for annealed and quenched random-graph
models see \cite{GGHP16,DGGHP16,KLS19,KLS20,KLS22}.

Our first result concerns high temperature, $0<\beta<1$.
To describe its deterministic centre and scale, recall the
classical Curie--Weiss mean-field equation. For an inverse
temperature $b\in(0,1)$ and external field $h\in\R$, let
$m(b,h)$ be the unique solution of the mean-field equation, and let
$\chi(b,h)$ be the corresponding susceptibility:
\begin{equation}\label{eq:intro-meanfield}
 m(b,h)=\tanh\bigl(bm(b,h)+h\bigr),\qquad
 \chi(b,h)=\frac{1-m(b,h)^2}{1-b\{1-m(b,h)^2\}}.
\end{equation}
For fixed $0<\beta<1$ and all sufficiently large $N$ in the dense regime,
$\beta_N<1$; write $m_N=m(\beta_N,h)$ and
$\chi_N=\chi(\beta_N,h)$.
The quantities $Nm_N$ and $N\chi_N$ are the deterministic
mean-field approximations to the mean and variance of $M_N$.

A Taylor estimate retaining the dependence on $p_N$ gives
\[
 |\beta_N-\beta|\le \frac{e\beta^3}{24(p_NN)^2}.
\]
Hence $p_NN\to\infty$ implies $\beta_N\to\beta<1$.
The homogeneous exchangeable-pair bound derived from
\cite[Theorem~1.4]{DE20} then applies uniformly with $b=\beta_N$
bounded away from $1$. Let $Z$ be a standard
normal random variable and $\Phi$ its distribution function.
We distinguish the mean-field and exact finite-volume standardisations:
\begin{equation}\label{eq:intro-standardizations}
 X_N=\frac{M_N-Nm_N}{\sqrt{N\chi_N}},\qquad
 Y_N=\frac{M_N-\E M_N}{\sqrt{\Var(M_N)}}.
\end{equation}
For these variables we obtain
\[
 d_{\mathrm K}(X_N,Z)\le \frac{C(\beta,h)}{\sqrt N},\qquad
 d_{\mathrm K}(Y_N,Z)\le \frac{C(\beta,h)}{\sqrt N}.
\]
The reduction also yields a sub-Gaussian concentration bound
for fixed $h$ via \cite{Cha07}, and at $h=0$ a bound with an
additional quartic tail term via \cite{CD10}.
The Cram\'er-type theorem of Chen, Fang and Shao \cite{CFS13}
yields, for the same $X_N$,
uniformly for $0\le x\le c_0N^{1/6}$ for some $c_0\in(0,1)$,
\[
 \frac{\P(X_N>x)}{1-\Phi(x)}
 =1+O_{\beta,h}\!\left(\frac{1+x^3}{\sqrt N}\right),
\]
and likewise for the lower tail.

Our second contribution concerns nominal criticality at zero external field $h=0$.  At $\beta=1$ a higher-order expansion gives, for $p_N\to0$,
\[
 \sqrt N(\beta_N-1)
 =\frac{1+o(1)}{24\sqrt{p_N^4N^3}}.
\]
The natural dilution threshold is therefore
\[
 p_N\asymp N^{-3/4}.
\]
If $p_N^4N^3\to c^{-2}\in(0,\infty)$, the limiting target has density
\[
 p_c(x)=Z_c^{-1}\exp\left\{\frac{c}{48}x^2-\frac{x^4}{12}\right\}.
\]
Here $Z_c=\int_{\mathbb R}\exp\{ct^2/48-t^4/12\}\,dt$.
The comparison densities $p_c$, $c\ge0$, have potential
and derivative
\[
 V_c(x)=\frac{x^4}{12}-\frac{c}{48}x^2,
 \qquad g_c(x)=V_c'(x)=\frac{x^3}{3}-\frac{c}{24}x.
\]
The function $g_c=-p_c'/p_c$ is called the \emph{Stein drift}.
For $z\in\R$ let $H_z(x)=\1_{\{x\le z\}}$ be the
indicator of the half-line $(-\infty,z]$. If $Y_c$ has density
$p_c$ and distribution function $F_c$, the corresponding Stein
equation is $f'(x)-g_c(x)f(x)=H_z(x)-F_c(z)$.
For $c=0$, the drift $g_0(x)=x^3/3$ is increasing and
has the sign of $x$. For $c>0$, however, $g_c(0)=0$ but
$g_c(x)<0$ on $(0,\sqrt{c/8})$. Thus the drift is neither
nondecreasing nor sign-preserving on the positive half-line;
correspondingly, $p_c$ has two modes rather than one.
These monotonicity and sign conditions underlie standard
Stein bounds for solutions of the indicator equation.
Since they fail here, the required bounds on the indicator Stein solutions and their derivatives
must be established directly for this double-well family.
Quartic critical laws and their quadratic perturbations are
classical \cite{EN78a,EN78b,ENR80}. Here the annealed
parameter satisfies $\sqrt N(\beta_N-1)\to c/24>0$,
selecting the nonmonotone branch. Its sign and coefficient
differ from the quenched crossover in \cite{KLS20}.
Our contribution is the quantitative Kolmogorov bound for this
annealed dilution window.
The exchangeable-pair approach to Stein's method is well established for Curie--Weiss and related spin systems.
The new difficulty here is to control indicator Stein solutions for a nonmonotone drift.

Shao, Zhang and Zhang \cite{SZZ24} develop normal and nonnormal
Berry--Esseen approximations for general Curie--Weiss models with
size-dependent inverse temperatures.  We do not claim that their framework
excludes the present parameter regime; our approach instead establishes
the explicit indicator Stein bounds needed for a Kolmogorov estimate
in the annealed double-well crossover.
Liu, Shao, Shi and Xu \cite{LSSX26} treat, in an edge--triangle
exponential random graph model, a quartic family with a quadratic
perturbation and obtain Wasserstein bounds.  Their Lipschitz-test
estimates do not directly provide the indicator-test estimates used here.

\begin{figure}[ht]
 \centering
 \includegraphics[width=0.82\textwidth]{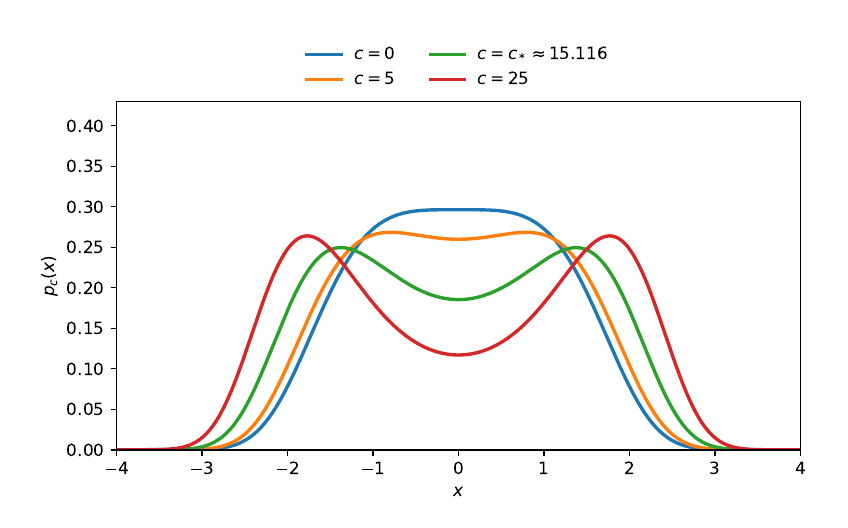}
 \caption{The crossover densities $p_c(x)=Z_c^{-1}\exp\{cx^2/48-x^4/12\}$ for representative values of $c$. Increasing $c$ turns the quartic critical law into an increasingly pronounced double-well density. The value $c_*\approx15.116$ marks the transition of the optimal uniform bound on the first derivatives of the indicator Stein solutions discussed in the appendix; it is not a threshold for bimodality.}
 \label{fig:crossover-densities}
\end{figure}

Figure~\ref{fig:crossover-densities} shows the passage from
the quartic law to a pronounced double well.

We solve the indicator Stein equation for $p_c$
explicitly. If $f_{c,z}$ denotes its bounded solution for
the indicator test $H_z$, we prove
\[
 \sup_z\|f_{c,z}\|_\infty=\frac{Z_c}{4},
 \qquad
 \sup_z\|f'_{c,z}\|_\infty
 \le 1+\frac{Z_c}{6}\left(\frac c{24}\right)^{3/2}.
\]
A target-specific exchangeable-pair argument then yields the quantitative crossover theorem
\[
 d_{\mathrm K}\left(\frac{M_N}{N^{3/4}},Y_c\right)
 \le C(c)\left\{N^{-1/2}
 +\left|\sqrt N(\beta_N-1)-\frac c{24}\right|\right\}.
\]
The deterministic mismatch term is essential: convergence $p_N^4N^3\to c^{-2}$ identifies the limiting law, but by itself does not imply an $N^{-1/2}$ rate for the effective parameter.

Section~\ref{sec:annealed-reduction} studies the effective
temperature, Section~\ref{sec:highT} the high-temperature
results, and Section~\ref{sec:critical} the critical window
and its Stein approximation. Appendix~\ref{app:optimal-K1}
analyzes the optimal uniform bound on the first derivatives of the indicator Stein solutions.

\section{The effective inverse temperature}
\label{sec:annealed-reduction}

The exact identity \eqref{eq:intro-betaN} allows a sharper estimate of the effective temperature. This is the key step behind the extension of the high-temperature Stein consequences to the full dense regime.  Its point is not a new approximation theorem, but the precise $p$-dependence of the Taylor remainder in the effective inverse temperature.  Although the high-temperature application below has $0<\beta<1$, the estimate itself is valid for every fixed $\beta>0$.

\begin{lemma}[Effective temperature]\label{lem:betaeff}
Let $\beta>0$, $0<p\le1$, and define $\beta_N$ by
\eqref{eq:intro-betaN} with $p_N=p$. If $pN\ge\beta$, then
\[
 \abs{\beta_N-\beta}
 \le \frac{e\beta^3}{24\,p^2N^2}
 =\frac{e\beta^3}{24\,(pN)^2}.
\]
The estimate is uniform in $p\in(0,1]$. Consequently, for every fixed $\beta>0$, the dense-regime condition $p_NN\to\infty$ implies $\beta_N\to\beta$.
If in addition $0<\beta<1$ and $\delta=(1-\beta)/2$, then for all sufficiently large $N$,
\[
 \beta_N\le 1-\delta=\frac{1+\beta}{2}<1.
\]
\end{lemma}

\begin{proof}
We keep track of the dependence on $p$ throughout the Taylor expansion.  This is the only
point at which the improvement over a $p$-independent cubic remainder enters.

Define
\[
 F_p(x)=\frac12\log\frac{1-p+pe^x}{1-p+pe^{-x}},
 \qquad x\in\R.
\]
By the definition \eqref{eq:intro-betaN} of the effective inverse temperature,
\begin{equation}\label{eq:betaN-via-Fp}
 \beta_N
 =2N F_p\!\left(\frac{\beta}{2pN}\right).
\end{equation}
Thus it suffices to control $F_p(x)-px$ for small $x$, uniformly in $p$.

For this purpose set
\[
 g_p(x)=\log(1-p+pe^x),
 \qquad
 r_p(x)=\frac{pe^x}{1-p+pe^x}.
\]
Then
\begin{equation}\label{eq:Fp-gp-relation}
 F_p(x)=\frac12\{g_p(x)-g_p(-x)\}.
\end{equation}
The function $r_p$ is the Bernoulli parameter obtained after the exponential tilt by $x$.
A direct differentiation gives
\begin{align}
 g_p'(x)&=r_p(x),\notag\\
 g_p''(x)&=r_p(x)\{1-r_p(x)\},\notag\\
 g_p'''(x)&=r_p(x)\{1-r_p(x)\}\{1-2r_p(x)\}.\label{eq:gp-third}
\end{align}
We now isolate the key uniform estimate.  If $|x|\le1/2$, then both $1$ and $e^x$ are at
least $e^{-1/2}$, and hence
\[
 1-p+pe^x=(1-p)\cdot1+p e^x\ge e^{-1/2}.
\]
At the same time $e^x\le e^{1/2}$.  Therefore
\[
 0\le r_p(x)
 =\frac{pe^x}{1-p+pe^x}
 \le \frac{p e^{1/2}}{e^{-1/2}}
 =ep,
 \qquad |x|\le\frac12.
\]
Since $0\le r_p(x)\le1$, formula \eqref{eq:gp-third} also gives
\[
 |g_p'''(x)|
 \le r_p(x)
 \le ep,
 \qquad |x|\le\frac12.
\]
The same bound of course holds with $x$ replaced by $-x$.

Differentiating \eqref{eq:Fp-gp-relation} three times yields
\[
 F_p'''(x)
 =\frac12\{g_p'''(x)+g_p'''(-x)\}.
\]
Combining the preceding two estimates therefore gives the
uniform estimate
\begin{equation}\label{eq:Fp-third-bound}
 |F_p'''(x)|\le ep,
 \qquad 0<p\le1,\quad |x|\le\frac12.
\end{equation}
This is the decisive point: the cubic derivative is of order $p$, not merely of order one.

For completeness, the lower derivatives at the origin are
\[
 F_p(0)=0,
 \qquad
 F_p'(0)=g_p'(0)=p,
 \qquad
 F_p''(0)=0,
\]
where the last identity also follows from the oddness of $F_p$.  Taylor's theorem with
remainder of order three now implies, for every $|x|\le1/2$,
\[
 \left|F_p(x)-px\right|
 \le \frac{|x|^3}{6}
      \sup_{|t|\le |x|}|F_p'''(t)|
 \le \frac{ep}{6}|x|^3.
\]

Finally, put $x_N=\beta/(2pN)$.
The assumption $pN\ge\beta$ is precisely what is needed to ensure
$0<x_N\le1/2$.  Hence the defining identity for $\beta_N$ and Taylor\'s bound give
\begin{align*}
 |\beta_N-\beta|
 &=2N|F_p(x_N)-px_N|
 \le \frac{Nep}{3}\left(\frac{\beta}{2pN}\right)^3\\
 &=\frac{e\beta^3}{24(pN)^2},
\end{align*}
which proves the claimed bound.

If now $p_NN\to\infty$, then the assumption $pN\ge\beta$ holds eventually and the
right-hand side of the bound tends to zero.  This proves the asserted convergence. Finally, $\beta_N\to\beta<1$ implies that, for
$\delta=(1-\beta)/2$, one has $|\beta_N-\beta|\le\delta$ for all sufficiently large $N$;
therefore $\beta_N\le\beta+\delta=(1+\beta)/2=1-\delta<1$, as claimed.
\end{proof}

\begin{remark}\label{rem:ADS-effective-temp}
The estimate in \cite{ADS26} gives
$\beta_N-\beta=O(\beta^3/(p^3N^2))$; it uses a cubic Taylor remainder
without retaining the factor $p$.  By contrast,
\eqref{eq:Fp-third-bound} yields $F_p(x)-px=O(px^3)$ uniformly for
$|x|\le 1/2$, and hence the sharper error $O(\beta^3/(pN)^2)$.
Thus $p_NN\to\infty$ suffices to control the effective inverse temperature.
This does not remove the stronger assumption $p_N^3N^2\to\infty$
from the other results of \cite{ADS26}, whose complex saddle-point
analysis requires more than real-parameter convergence of $\beta_N$.
\end{remark}

\section{High-temperature normal approximation and further Stein consequences}\label{sec:highT}
\subsection{The classical Curie--Weiss Berry--Esseen input}
\label{subsec:classical-CW-input}

By the exact annealed reduction, the high-temperature problem reduces to
a classical Curie--Weiss model with inverse temperature $b=\beta_N$.
We establish the required compact-uniform Berry--Esseen bound directly
in the homogeneous model, using the one-dimensional exchangeable-pair
inequality \cite[Theorem~1.4]{DE20}.

We use the mean-field quantities $m(b,u)$ and $\chi(b,u)$ defined
in \eqref{eq:intro-meanfield}, with $h$ there replaced by $u$.
In the homogeneous case, these agree with the deterministic centre and
susceptibility in \cite[Theorem~1.3]{DE20}; they must not be confused
with the random total magnetisation $M_N$. The law $CW(b,u)$ is the
classical Gibbs law defined after \eqref{eq:intro-exact-CW}.

\begin{lemma}[Compact-uniform classical Berry--Esseen bound]\label{lem:DE-uniform}
Let $0<b_-<b_+<1$ and let $H\subset\R$ be compact.  There is a finite constant
$C=C(b_-,b_+,H)$ such that, for every $b\in[b_-,b_+]$, $u\in H$ and every $N$,
\begin{equation}\label{eq:DE-classical-BE}
 d_K\left(
 \frac{M_N-Nm(b,u)}{\sqrt{N\chi(b,u)}},Z
 \right)\le \frac{C}{\sqrt N}
\end{equation}
under $CW(b,u)$ as defined after \eqref{eq:intro-exact-CW}.
\end{lemma}

\begin{proof}
We verify uniformity directly in the homogeneous model, using the
one-dimensional case of \cite[Theorem~1.4]{DE20}.  Constants below
may depend on $b_-,b_+$ and $H$, but not on $b,u,N$.
Put $m=m(b,u)$, $q=1-m^2$, $a=1-bq$ and
$X=(M_N-Nm)/\sqrt{N\chi}$, where $\chi=q/a$.
Compactness gives $q\ge q_->0$, $a\ge1-b_+>0$ and
$0<\chi_-\le\chi\le\chi_+<\infty$.

Resample a uniformly selected spin from its conditional Gibbs law,
obtaining $M_N'$ and the exchangeable pair
$X'=(M_N'-Nm)/\sqrt{N\chi}$; write $D=X-X'$ and
$s=M_N/N$.  Conditional on the configuration,
\[
 \E(D\mid\sigma)=\frac{1}{\sqrt{N\chi}}
 \left(s-\frac1N\sum_{i=1}^N
 \tanh\bigl(u+b(s-\sigma_i/N)\bigr)\right),
\]
while
$\E(D^2\mid\sigma)=2(N\chi)^{-1}
\{1-N^{-1}\sum_i\sigma_i\tanh(u+b(s-\sigma_i/N))\}$.
Since $\tanh$ is $1$-Lipschitz, replacing its arguments by
$u+bs$ incurs an error at most $b_+/N$.

First, $\sup_{b,u,N}\E X^2<\infty$.  Indeed,
$f(s)=s-\tanh(u+bs)$ satisfies $f(m)=0$ and
$f'(s)\ge1-b_+$.  Hence the conditional drift formula gives
$X\E(D\mid\sigma)\ge(1-b_+)X^2/N-C|X|N^{-3/2}$.
Exchangeability yields $\E\{X\E(D\mid\sigma)\}
=\E D^2/2\le2/(N\chi_-)$; Young's inequality proves
$\E X^2\le C$ uniformly.

Taylor expansion at $m$, using $m=\tanh(u+bm)$ and the
bounded second derivative of $\tanh$, now gives
\[
 \E(D\mid X)=\frac aN(X+R),
 \qquad |R|\le\frac C{\sqrt N}(1+X^2).
\]
Moreover, conditioning the second-moment identity on $X$ gives
\[
 \E(D^2\mid X)=\frac{2a}{N}
 +O\bigl(N^{-3/2}(1+|X|)\bigr),
\]
with a uniform pointwise error.
Finally, since $D\in\{0,\pm2/\sqrt{N\chi}\}$,
$|D|D=2D/\sqrt{N\chi}$.  Consequently, with
$\lambda=a/N$,
\[
 \E\left|1-\frac{\E(D^2\mid X)}{2\lambda}\right|
 +\frac1\lambda\E\big|\E(|D|D\mid X)\big|
 +\E|R|\le\frac C{\sqrt N}.
\]
The one-dimensional bound \cite[Theorem~1.4, equation~(1.13)]{DE20}
proves \eqref{eq:DE-classical-BE}.
\end{proof}

For the comparison with Apostel--D\"oring--Schubert we also need to pass from the explicit
mean-field normalization in \eqref{eq:DE-classical-BE} to the exact finite-volume mean and
variance.  The following standard finite-volume estimate is included to make that passage
explicit.

\begin{lemma}[Uniform finite-volume matching]\label{lem:CW-finite-matching}
Let $0<b_-<b_+<1$ and let $H\subset\R$ be compact.  There is a finite constant
$C=C(b_-,b_+,H)$ such that, uniformly for $b\in[b_-,b_+]$, $u\in H$,
\begin{equation}\label{eq:CW-mean-match}
 \left|\E_{b,u}M_N-Nm(b,u)\right|\le C,
\end{equation}
and
\begin{equation}\label{eq:CW-var-match}
 \left|\Var_{b,u}(M_N)-N\chi(b,u)\right|\le C.
\end{equation}
\end{lemma}

\begin{proof}
We avoid differentiating a parameter-dependent Laplace remainder.  Introduce the
Hubbard--Stratonovich auxiliary variable $X$ through the joint law on
$\{-1,1\}^N\times\R$ with density proportional to
\[
 \sqrt{\frac{N}{2\pi b}}
 \exp\left\{-\frac{N}{2b}x^2+(x+u)M_N(\sigma)\right\}\,dx.
\]
Integrating out $x$ gives exactly the Curie--Weiss weight
$\exp\{bM_N^2/(2N)+uM_N\}$.  Conversely, conditional on $M_N$,
completion of the square gives
\[
 X\mid M_N
 \sim \mathcal N\!\left(\frac{bM_N}{N},\frac bN\right).
\]
Hence, if $\nu_{N,b,u}$ denotes the marginal law of $X$, then
\begin{equation}\label{eq:CW-HS-moment-identities}
 \E_{\nu}X=\frac bN\E_{b,u}M_N,
 \qquad
 \Var_{\nu}(X)=\frac bN+\frac{b^2}{N^2}\Var_{b,u}(M_N).
\end{equation}
The $X$-marginal has density proportional to
\[
 \exp\{N\phi_{b,u}(x)\},
 \qquad
 \phi_{b,u}(x)=-\frac{x^2}{2b}+\log\bigl(2\cosh(x+u)\bigr).
\]

For $b<1$, $\phi_{b,u}$ has the unique maximizer
$x_*(b,u)=bm(b,u)$.  Put
\[
 q(b,u)=1-m(b,u)^2,
 \qquad
 \kappa(b,u):=-\phi_{b,u}''(x_*)
 =\frac{1-bq(b,u)}{b}.
\]
The Laplace estimates needed below are uniform on
$[b_-,b_+]\times H$.  Indeed,
\[
 \phi_{b,u}''(x)
 =-\frac1b+\operatorname{sech}^2(x+u)
 \le -\frac{1-b_+}{b_+}=: -\kappa_0<0,
\]
so $\phi_{b,u}$ is globally and uniformly strictly concave, and therefore
\[
 \phi_{b,u}(x)
 \le \phi_{b,u}(x_*)-\frac{\kappa_0}{2}(x-x_*)^2.
\]
Moreover, $\kappa$ is bounded above and below by positive constants, while all
derivatives of $\phi_{b,u}$ of order at least three are globally bounded,
uniformly in $(b,u)$.

With the change of variables $y=\sqrt N\,(x-x_*)$, for $r=0,1,2$ write
\[
 I_{r,N}
 =\int_{\R} y^r
 \exp\left\{N\bigl(\phi_{b,u}(x_*+y/\sqrt N)-\phi_{b,u}(x_*)\bigr)\right\}\,dy.
\]
We now record the Laplace expansion to the order needed for the variance estimate.  Since
$\phi'_{b,u}(x_*)=0$ and $\phi''_{b,u}(x_*)=-\kappa$, Taylor's theorem gives, uniformly in
$(b,u)\in[b_-,b_+]\times H$,
\begin{align*}
 &N\bigl(\phi_{b,u}(x_*+y/\sqrt N)-\phi_{b,u}(x_*)\bigr)\\
 &\qquad=-\frac{\kappa y^2}{2}
 +\frac{\phi_{b,u}'''(x_*)}{6\sqrt N}y^3
 +\frac{\phi_{b,u}^{(4)}(x_*)}{24N}y^4
 +O\left(\frac{|y|^5}{N^{3/2}}\right)
\end{align*}
whenever $|y|\le N^{1/10}$.  All constants are uniform because the derivatives of
$\phi_{b,u}$ of order at least three are uniformly bounded.  On this central region the
last three terms in the Taylor expansion above are uniformly small, and hence
\begin{align*}
 &\exp\left\{N\bigl(\phi_{b,u}(x_*+y/\sqrt N)-\phi_{b,u}(x_*)\bigr)\right\}\\
 &\qquad=e^{-\kappa y^2/2}
 \left[
 1+\frac{\phi_{b,u}'''(x_*)}{6\sqrt N}y^3
 +O\left(\frac{y^4+y^6}{N}\right)
 \right].
\end{align*}
The error in this expansion is uniform.  On the complement
$|y|>N^{1/10}$, the global concavity bound gives
Gaussian domination by $e^{-\kappa_0y^2/2}$, so the corresponding contribution to each
$I_{r,N}$, $r=0,1,2$, is smaller than any fixed power of $N^{-1}$.

Integrating this expansion over the symmetric central interval now
shows exactly where the improved orders come from.  The possible $N^{-1/2}$ correction to
$I_{0,N}$ is proportional to
\[
 \int y^3e^{-\kappa y^2/2}\,dy,
\]
and the corresponding correction to $I_{2,N}$ is proportional to
\[
 \int y^5e^{-\kappa y^2/2}\,dy.
\]
Both vanish by oddness.  For $I_{1,N}$, in contrast, the first nonzero correction is of order
$N^{-1/2}$ and involves $\int y^4e^{-\kappa y^2/2}\,dy$.  Consequently, uniformly on the
parameter compact set,
\begin{align*}
 I_{0,N}
 &=\int_{\R}e^{-\kappa y^2/2}\,dy+O(N^{-1}),\\
 I_{1,N}&=O(N^{-1/2}),\\
 I_{2,N}
 &=\int_{\R}y^2e^{-\kappa y^2/2}\,dy+O(N^{-1}).
\end{align*}

The moments of $X$ are ratios of these three integrals.  More precisely,
\[
 \E_{\nu}X
 =x_*+\frac1{\sqrt N}\frac{I_{1,N}}{I_{0,N}},
\]
and
\[
 \Var_{\nu}(X)
 =\frac1N\left\{
 \frac{I_{2,N}}{I_{0,N}}
 -\left(\frac{I_{1,N}}{I_{0,N}}\right)^2
 \right\}.
\]
Since $\kappa$ is bounded above and below by positive constants, the Gaussian normalizing
integral $I_{0,N}$ is uniformly bounded above and away from zero.  Furthermore,
\[
 \frac{\int_{\R}y^2e^{-\kappa y^2/2}\,dy}
      {\int_{\R}e^{-\kappa y^2/2}\,dy}
 =\frac1\kappa.
\]
These estimates imply
\[
 \E_{\nu}X=x_*+O(N^{-1}),
 \qquad
 \Var_{\nu}(X)=\frac1{N\kappa}+O(N^{-2}),
\]
with constants uniform in $(b,u)$.  Notice that the cancellation of the odd
$N^{-1/2}$ term in $I_{0,N}$ and $I_{2,N}$ is precisely what yields the $O(N^{-2})$
remainder in the variance, which is needed after the $N^2$ rescaling below.  Using $x_*=bm$ in the first identity of
\eqref{eq:CW-HS-moment-identities} gives
\[
 \E_{b,u}M_N
 =\frac Nb\E_{\nu}X
 =Nm(b,u)+O(1),
\]
uniformly because $b\ge b_->0$.

Finally, the second identity in \eqref{eq:CW-HS-moment-identities}, together with
$\kappa=(1-bq)/b$, immediately gives
\[
 \Var_{b,u}(M_N)=\frac{N}{b^2}\left(\frac1\kappa-b\right)+O(1)
 =N\chi(b,u)+O(1).
\]
The bounds are uniform for all sufficiently large $N$; enlarging $C$ covers the finitely many remaining values.  This proves \eqref{eq:CW-mean-match}--\eqref{eq:CW-var-match}.
\end{proof}

\subsection{The full dense annealed theorem}

\begin{theorem}[Berry--Esseen in the full dense annealed regime]\label{thm:main}
Let $0<\beta<1$ and $h\in\R$ be fixed.  Let $p_N\in(0,1]$ satisfy
\[
 p_NN\longrightarrow\infty.
\]
Let $M_N$ have the annealed dilute Curie--Weiss distribution
\eqref{eq:intro-exact-CW}. For $X_N,Y_N$ as in
\eqref{eq:intro-standardizations}, there exist $N_0$ and
$C=C(\beta,h)<\infty$ such that, for $N\ge N_0$,
\begin{equation}\label{eq:main-BE-mf}
 d_K(X_N,Z)\le\frac{C}{\sqrt N},
\end{equation}
and, with exact finite-volume standardization,
\begin{equation}\label{eq:main-BE}
 d_K(Y_N,Z)\le\frac{C}{\sqrt N}.
\end{equation}
\end{theorem}

\begin{proof}
By \eqref{eq:intro-exact-CW}, the annealed dilute model is exactly the classical
Curie--Weiss model with parameters $(\beta_N,h)$.  Lemma~\ref{lem:betaeff} gives
\[
 |\beta_N-\beta|\le \frac{e\beta^3}{24(p_NN)^2},
\]
so $p_NN\to\infty$ implies $\beta_N\to\beta<1$.  Choose
$0<b_-<\beta<b_+<1$ so that $\beta_N\in[b_-,b_+]$ for all sufficiently large $N$, and
apply Lemma~\ref{lem:DE-uniform} with $H=\{h\}$.  This proves \eqref{eq:main-BE-mf}.

Set
\[
 a_N=\E X_N,
 \qquad
 s_N^2=\Var(X_N).
\]
Since $\chi(b,h)$ is continuous and strictly positive on $[b_-,b_+]$, Lemma~\ref{lem:CW-finite-matching}
gives
\[
 |a_N|\le C N^{-1/2},
 \qquad
 |s_N^2-1|\le C N^{-1},
\]
and hence $|s_N-1|\le C N^{-1}$ for large $N$.  Moreover,
\[
 Y_N=\frac{X_N-a_N}{s_N}.
\]
Therefore, for every $z$,
\begin{align*}
 \left|\P\left(\frac{X_N-a_N}{s_N}\le z\right)-\Phi(z)\right|
 &\le d_K(X_N,Z)+|\Phi(a_N+s_Nz)-\Phi(z)|.
\end{align*}
For $1/2\le s\le2$,
\[
 \sup_z|\Phi(a+sz)-\Phi(z)|\le C\{|a|+|s-1|\}.
\]
Here $\varphi=\Phi'$ denotes the standard normal density.
The bound follows from the mean-value theorem and
$\sup_x|x|\varphi(x)<\infty$. Combining it with
\eqref{eq:main-BE-mf} proves \eqref{eq:main-BE}.
\end{proof}

\subsection{Two further high-temperature consequences}
\label{subsec:highT-further}

The exact Curie--Weiss representation \eqref{eq:intro-exact-CW} also allows us to obtain
concentration bounds and a classical Cram\'er-type normal-tail approximation.
We record two concentration inequalities: the first, based on \cite{Cha07}, gives a
sub-Gaussian bound for fixed external field; the second, based on \cite{CD10},
applies at zero field and contains an additional quartic term in the exponent.

\begin{proposition}[Concentration in the full dense regime]\label{prop:highT-concentration}
Under the assumptions of Theorem~\ref{thm:main}, put
$\overline M_N=M_N/N$ and recall $m_N=m(\beta_N,h)$.
Then, for all sufficiently large $N$ and every $t\ge0$,
\begin{equation}\label{eq:highT-concentration-mf}
 \P\left(
 |\overline M_N-m_N|
 \ge \frac{\beta_N}{N(1-\beta_N)}
      +\frac{t}{(1-\beta_N)\sqrt N}
 \right)
 \le 2\exp\left\{-\frac{t^2}{4(1+\beta_N)}\right\}.
\end{equation}
Consequently there are constants $c,C\in(0,\infty)$, depending only on $(\beta,h)$, such that
for $Y_N$ from \eqref{eq:intro-standardizations} one has
\begin{equation}\label{eq:highT-concentration-exact}
 \P(|Y_N|\ge x)
 \le 2\exp\{-c(x-CN^{-1/2})_+^2\},
 \qquad x\ge0,
\end{equation}
where $y_+=\max\{y,0\}$.
\end{proposition}

\begin{proof}
By the exact reduction, it is enough to work under $CW(b,h)$ with $b=\beta_N<1$.
Construct the heat-bath exchangeable pair $(\sigma,\sigma')$ by choosing an index
$I$ uniformly from $[N]$ and resampling $\sigma_I$ from its conditional Gibbs law.
For this pair put
\[
 F(\sigma,\sigma')=M_N(\sigma)-M_N(\sigma'),\qquad
 f(\sigma)=\E\bigl[F(\sigma,\sigma')\mid\sigma\bigr].
\]
Then $F$ is antisymmetric and
\[
 f(\sigma)
 =\frac1N\sum_{i=1}^N
 \left\{\sigma_i-\tanh\left(h+b\frac{M_N-\sigma_i}{N}\right)\right\}.
\]
If two configurations differ at one site, then the explicit sum and the
$1$-Lipschitz property of $\tanh$ give
\[
 |f(\sigma)-f(\widetilde\sigma)|\le \frac{2(1+b)}{N}.
\]
Since $|F|\le2$, the conditional variance proxy in Chatterjee's exchangeable-pair
concentration theorem \cite[Theorem~1.5]{Cha07} satisfies
\[
 \frac12\E\left(
 |(f(\sigma)-f(\sigma'))F(\sigma,\sigma')|
 \mid\sigma\right)
 \le \frac{2(1+b)}{N}.
\]
Hence
\[
 \P(|f(\sigma)|\ge s)
 \le 2\exp\left\{-\frac{Ns^2}{4(1+b)}\right\},
 \qquad s\ge0.
\]
On the other hand,
\[
 \left|
 f(\sigma)-\{\overline M_N-\tanh(h+b\overline M_N)\}
 \right|\le \frac bN.
\]
The function $T(x)=\tanh(h+bx)$ is $b$-Lipschitz, and its fixed point
is $m(b,h)=m_N$. Hence
\[
 (1-b)|\overline M_N-m_N|
 \le |\overline M_N-T(\overline M_N)|.
\]
Combining these bounds with the preceding concentration estimate proves
\eqref{eq:highT-concentration-mf}.

Finally, Lemma~\ref{lem:CW-finite-matching} gives
$|\E M_N-Nm_N|=O(1)$ and $\Var(M_N)\asymp N$ uniformly for large $N$.
Since $\beta_N$ stays in a compact subinterval of $(0,1)$, the explicit bound
\eqref{eq:highT-concentration-mf} therefore implies \eqref{eq:highT-concentration-exact}
after changing the constants $c,C$.
\end{proof}

\begin{corollary}[Chatterjee--Dey Gaussian--quartic concentration at zero field]\label{prop:highT-concentration-CD}
Under the assumptions of Theorem~\ref{thm:main}, assume in addition that $h=0$.
Then, for all sufficiently large $N$ and every $t\ge0$,
\[
 \P\left(
 3(1-\beta_N)\overline M_N^2+\beta_N^3\overline M_N^4\ge t
 \right)
 \le 2\exp\left\{-\frac{Nt}{160}\right\}.
\]
Equivalently, with the mean-field high-temperature normalization
\[
 \widetilde M_N
 :=\sqrt{N(1-\beta_N)}\,\overline M_N
 =\sqrt{\frac{1-\beta_N}{N}}\,M_N,
\]
one has, for every $x\ge0$,
\[
 \P(|\widetilde M_N|\ge x)
 \le
 2\exp\left\{
 -\frac{3x^2}{160}
 -\frac{\beta_N^3x^4}{160N(1-\beta_N)^2}
 \right\}.
\]
This is \cite[Proposition~5]{CD10}, applied after the exact annealed reduction with $b=\beta_N$.
\end{corollary}

\begin{remark}\label{rem:CD-concentration}
The two concentration statements are complementary.  Proposition~\ref{prop:highT-concentration}
allows arbitrary fixed $h$ and transfers to exact centering and variance standardization, whereas
Proposition~\ref{prop:highT-concentration-CD} is a zero-field finite-$N$ inequality whose quartic
term retains the change of tail geometry as the classical critical point is approached.  The latter
result assumes $b\le1$ and therefore does not cover the annealed positive-tilt crossover, where
$\beta_N>1$ for all sufficiently large $N$ when $p_N\to0$.
\end{remark}

\begin{theorem}[Cram\'er-type normal tails in the full dense regime]
\label{thm:highT-cramer}
Under the assumptions of Theorem~\ref{thm:main}, let $X_N$ and $Y_N$ be as in
\eqref{eq:intro-standardizations}. There are $N_0$ and $C=C(\beta,h)<\infty$ such that, for $N\ge N_0$ and
$0\le x\le N^{1/6}$,
\begin{equation}\label{eq:highT-cramer-mf}
 \left|
 \frac{\P(X_N>x)}{1-\Phi(x)}-1
 \right|
 \le \frac{C(1+x^3)}{\sqrt N}.
\end{equation}
The analogous estimate holds for the lower tail. Moreover, for the exactly
standardised $Y_N$ there are $c_0\in(0,1)$ and $C<\infty$ such that
\begin{equation}\label{eq:highT-cramer-exact}
 \left|
 \frac{\P(Y_N>x)}{1-\Phi(x)}-1
 \right|
 \le \frac{C(1+x^3)}{\sqrt N},
 \qquad 0\le x\le c_0N^{1/6},
\end{equation}
and likewise for the lower tail.
\end{theorem}

\begin{proof}
By \eqref{eq:intro-exact-CW}, it suffices to apply
\cite[Proposition~4.3]{CFS13} to $CW(b,h)$ with $b=\beta_N$.
The field parameter in that paper is multiplied by $b$, so its
$h_{\mathrm{CFS}}$ equals $h/b$. Its mean-field equation and variance
normalisation are then precisely \eqref{eq:intro-meanfield} and
\eqref{eq:intro-standardizations}. Since $\beta_N\to\beta\in(0,1)$, $b$ lies eventually in a
fixed compact interval $[b_-,b_+]\subset(0,1)$, and the field
$h/b$ used in \cite{CFS13} stays bounded.  Here are the uniformity
details for \cite[Section~6.3]{CFS13}.  The magnetisation potential
$G_{b,h}(s)=I(s)-bs^2/2-hs$ satisfies
$G_{b,h}''(s)=(1-s^2)^{-1}-b\ge1-b_+$ for $|s|<1$;
its unique minimum $m(b,h)$ remains in a compact subinterval of
$(-1,1)$.  Thus a neighbourhood of this minimum can be chosen
uniformly, with a strictly positive uniform potential gap outside it.
The binomial weights then give the exponentially small truncation
error in that proof with a uniform exponent.  Inside the neighbourhood,
its regression coefficient is $\lambda=\{1-b(1-m^2)\}/N$, hence
$N\lambda\ge1-b_+$.  The normalising susceptibility
$\chi=(1-m^2)/\{1-b(1-m^2)\}$ is bounded above and away from zero,
and the first three derivatives of $\tanh$ are bounded.  The explicit
conditional-variance and Taylor-remainder estimates in
\cite[Section~6.3]{CFS13} therefore have constants bounded uniformly
on this parameter set; the truncation neighbourhood can also be
chosen so that the condition $\delta_2|W|\le1/2$ used there holds
with the same constants.  Consequently the constant in
\cite[Proposition~4.3]{CFS13} is uniform in our parameter range. We obtain
\eqref{eq:highT-cramer-mf} first with $\P(X_N\ge x)$ in place of
$\P(X_N>x)$, exactly as stated in \cite{CFS13}.

To pass to strict tails, note that $X_N$ is supported on a lattice
of mesh $d_N=2/\sqrt{N\chi_N}=O(N^{-1/2})$. For $x\ge0$,
\[
 \P(X_N\ge x+d_N)\le\P(X_N>x)\le\P(X_N\ge x).
\]
The normal tail ratio
$(1-\Phi(x+d_N))/(1-\Phi(x))=1+O((1+x)d_N)$
is uniform for $0\le x\le N^{1/6}$. Applying the cited bound
at $x$ and $x+d_N$ gives \eqref{eq:highT-cramer-mf} whenever
$x+d_N\le N^{1/6}$. In the remaining interval
$N^{1/6}-d_N<x\le N^{1/6}$, the relative error allowed by
\eqref{eq:highT-cramer-mf} is bounded below by a positive constant;
the upper bound from $\P(X_N>x)\le\P(X_N\ge x)$ and the trivial
lower bound $\P(X_N>x)\ge0$ therefore suffice after enlarging $C$.
This proves \eqref{eq:highT-cramer-mf}. Spin reversal, which changes $h$ to $-h$,
gives the lower tail.

For exact finite-volume standardisation, let $a_N=\E X_N$ and
$s_N^2=\Var(X_N)$. Lemma~\ref{lem:CW-finite-matching} gives
\[
 Y_N=\frac{X_N-a_N}{s_N},\qquad
 |a_N|=O(N^{-1/2}),\quad |s_N-1|=O(N^{-1}).
\]
Thus $\P(Y_N>x)=\P(X_N>a_N+s_Nx)$.
Choose $c_0\in(0,1)$ so that, for all large $N$,
$a_N+s_Nx\le N^{1/6}$ whenever $0\le x\le c_0N^{1/6}$.
Whenever $a_N+s_Nx\ge0$, the mean-field tail estimate applies,
and the normal Mills ratio yields, uniformly in this range,
\[
 \frac{1-\Phi(a_N+s_Nx)}{1-\Phi(x)}
 =1+O\!\left(\frac{1+x^3}{\sqrt N}\right).
\]
If $a_N+s_Nx<0$, necessarily $x=O(N^{-1/2})$. In this small
range the Berry--Esseen bound \eqref{eq:main-BE-mf}, together with
the Lipschitz continuity of $\Phi$ and $1-\Phi(x)\asymp1$,
gives the same relative error. This proves
\eqref{eq:highT-cramer-exact}. The lower-tail argument is analogous.
\end{proof}

\section{Nominal criticality and the annealed crossover}
\label{sec:critical}

The critical temperature deserves a separate treatment.  The results of \cite{ADS26} are stated for $0<\beta<1$; their Statulevi\v{c}ius theorem does not include $\beta=1$.  Throughout the critical analysis we assume zero external field $h=0$.  Because the annealed measure is then exactly $CW(\beta_N,0)$, the critical behavior at nominal $\beta=1$ can be read from the asymptotics of the effective inverse temperature $\beta_N$.

At criticality the relevant Curie--Weiss window is
\[
 \sqrt N\,(\beta_N-1)=O(1).
\]
To see why, write $m=M_N/N$. Under $CW(\beta_N,0)$, a magnetisation
value $M_N=Nm$ has weight proportional to
$\binom{N}{(N+Nm)/2}\exp\{N\beta_Nm^2/2\}$.
The binomial factor contributes the entropy term $-NI(m)$ to the
leading exponent, where
\[
 I(m)=\frac{1+m}{2}\log(1+m)+\frac{1-m}{2}\log(1-m).
\]
Consequently, apart from terms independent of $m$ and the Stirling
prefactor, the exponent is $N\{\beta_Nm^2/2-I(m)\}$.
Since $I(m)=m^2/2+m^4/12+O(m^6)$ as $m\to0$, setting
$m=N^{-1/4}x$ for fixed $x$ gives
\[
 N\left\{\frac{\beta_Nm^2}{2}-I(m)\right\}
 =\frac{\sqrt N(\beta_N-1)}2x^2-\frac{x^4}{12}+O(N^{-1/2}).
\]
Thus $\sqrt N(\beta_N-1)$ determines whether a quadratic perturbation
remains visible on the $N^{3/4}$ magnetisation scale.

\subsection{Critical expansion of the effective temperature}

Lemma~\ref{lem:betaeff}, now applied at $\beta=1$, already shows under $p_NN\to\infty$ that $\beta_N\to1$.  For the critical discussion, however, this absolute-error estimate is not sufficient: the sign and the leading coefficient of $\beta_N-1$ are essential.

\begin{lemma}[Critical expansion of the effective inverse temperature]\label{lem:critical-betaeff}
Let the nominal inverse temperature be $\beta=1$, and assume $p_NN\to\infty$. Then
\begin{equation}\label{eq:critical-betaeff-expansion}
 \beta_N
 =1+\frac{(1-p_N)(1-2p_N)}{24p_N^2N^2}
 +O\left(\frac1{p_N^4N^4}\right).
\end{equation}
In particular, if $p_N\to0$,
\begin{equation}\label{eq:critical-gamma}
 \theta_N:=\sqrt N(\beta_N-1)
 =\frac{1+o(1)}{24p_N^2N^{3/2}}
 =\frac{1+o(1)}{24\sqrt{p_N^4N^3}}.
\end{equation}
\end{lemma}

\begin{proof}
We use the notation $F_p$, $g_p$ and $r_p$ from the proof of
Lemma~\ref{lem:betaeff}, taking $p=p_N$ when applying the expansion below.
In particular, $F_p(x)=\{g_p(x)-g_p(-x)\}/2$ is odd,
$F_p'(0)=p$, and the already established formula
\eqref{eq:gp-third}, together with $r_p(0)=p$, gives
$F_p'''(0)=g_p'''(0)=p(1-p)(1-2p)$.
The third derivative at the origin is therefore known exactly;
to obtain a uniform fifth-order remainder we need two further derivatives.
They give
\[
 g_p^{(5)}(x)
 =r_p(x)\{1-r_p(x)\}\{1-2r_p(x)\}
   \{12r_p(x)^2-12r_p(x)+1\}.
\]
Fix $x_0>0$. As in Lemma~\ref{lem:betaeff}, the estimate
$r_p(x)\le e^{2x_0}p$ holds whenever $|x|\le x_0$ and
$0<p\le1$. Since $0\le r_p(x)\le1$, the remaining factors
are bounded uniformly. Thus $|g_p^{(5)}(x)|\le C_{x_0}p$;
the same bound holds at $-x$ and hence for
$F_p^{(5)}(x)=\{g_p^{(5)}(x)+g_p^{(5)}(-x)\}/2$.
By oddness and Taylor's theorem,
\[
 F_p(x)=px+\frac{p(1-p)(1-2p)}6x^3+R_p(x),
 \qquad |R_p(x)|\le C_{x_0}p|x|^5
\]
for $|x|\le x_0$, with a constant independent of $p$.

At $\beta=1$, identity \eqref{eq:betaN-via-Fp} reads
$\beta_N=2NF_{p_N}(x_N)$ with $x_N=1/(2p_NN)$.
Since $p_NN\to\infty$, we have $x_N\to0$ and may apply the
expansion for all sufficiently large $N$. Its linear term is
$2Np_Nx_N=1$, its cubic term is
$(1-p_N)(1-2p_N)/(24p_N^2N^2)$, and its remainder satisfies
$2N|R_{p_N}(x_N)|\le C/(p_N^4N^4)$.
This proves \eqref{eq:critical-betaeff-expansion}.

Multiplying by $\sqrt N$ gives
\[
 \theta_N=\frac{(1-p_N)(1-2p_N)}{24p_N^2N^{3/2}}
       +O\left(\frac1{p_N^4N^{7/2}}\right).
\]
The remainder divided by $1/(p_N^2N^{3/2})$ is
$O((p_NN)^{-2})=o(1)$. If also $p_N\to0$, then
$(1-p_N)(1-2p_N)=1+o(1)$, proving \eqref{eq:critical-gamma}.
\end{proof}

\begin{remark}
When $p_N\to0$, the leading correction in \eqref{eq:critical-betaeff-expansion} is \emph{positive}.  Thus the annealed model at nominal $\beta=1$ is effectively a slightly \emph{supercritical} classical Curie--Weiss model.  This sign will be responsible for the positive quadratic tilt in the crossover distribution below.
\end{remark}

\subsection{The three effective-temperature regimes}

Write $r_N:=p_N^4N^3=(p_NN^{3/4})^4$.  The following three
regimes follow from Lemma~\ref{lem:critical-betaeff}.

\begin{proposition}[Effective-temperature regimes at nominal $\beta=1$]\label{prop:annealed-critical-regimes}
Assume $h=0$ and $p_NN\to\infty$, and recall
$\theta_N=\sqrt N(\beta_N-1)$ from \eqref{eq:critical-gamma}.
\begin{enumerate}
\item If $r_N\to\infty$, then $\theta_N\to0$.
\item If $r_N\to c^{-2}$ for some $c\in(0,\infty)$, then $p_N\to0$
and $\theta_N\to c/24$.
\item If $r_N\to0$, then $p_N\to0$, $\theta_N\to\infty$, and
$\beta_N>1$ for all sufficiently large $N$, with $\beta_N\to1$.  The nonzero solutions $\pm m_N^*$ of
$m=\tanh(\beta_Nm)$ satisfy
\[
 m_N^*\sim\sqrt{3(\beta_N-1)}\sim\frac{1}{\sqrt8\,p_NN},
 \qquad N^{1/4}m_N^*\sim\frac{1}{\sqrt8\,p_NN^{3/4}}\to\infty.
\]
Their local susceptibility and the associated fluctuation scale satisfy
\[
 \chi_N^+:=\frac{1-(m_N^*)^2}{1-\beta_N[1-(m_N^*)^2]}
 \sim\frac{1}{2(\beta_N-1)}\sim12p_N^2N^2,
 \qquad \sqrt{N\chi_N^+}\sim\sqrt{12}\,p_NN^{3/2}.
\]
These are deterministic equilibrium and susceptibility asymptotics;
no conditional central limit theorem is claimed.
\end{enumerate}
\end{proposition}

\begin{proof}
By \eqref{eq:critical-betaeff-expansion},
\begin{equation}\label{eq:theta-full-expansion}
 \theta_N
 =\frac{(1-p_N)(1-2p_N)}{24p_N^2N^{3/2}}
 +O\!\left(\frac1{p_N^4N^{7/2}}\right).
\end{equation}
Since $0<p_N\le1$, this gives
$|\theta_N|\le C\{r_N^{-1/2}+(r_N\sqrt N)^{-1}\}\to0$
when $r_N\to\infty$, without requiring $p_N\to0$.
If $r_N\to c^{-2}\in(0,\infty)$ or $r_N\to0$, then
$p_N\to0$; \eqref{eq:critical-gamma} gives respectively
$\theta_N\to c/24$ and $\theta_N\to\infty$.
In the latter case \eqref{eq:critical-betaeff-expansion} also yields
$\beta_N-1\sim(24p_N^2N^2)^{-1}>0$ and $\beta_N\to1$.

For $b>1$, the equation $m=\tanh(bm)$ has exactly one positive
solution: equivalently, $b=\operatorname{arctanh}(m)/m$,
whose right-hand side increases strictly from $1$ to infinity on
$(0,1)$.  Thus $m_N^*\to0$, and the expansion
$\operatorname{arctanh}(m)/m=1+m^2/3+O(m^4)$ gives
$(m_N^*)^2\sim3(\beta_N-1)$.  Substituting the preceding
asymptotic for $\beta_N-1$ gives both equilibrium-point scales.
Finally, the denominator of $\chi_N^+$ is
$1-\beta_N+\beta_N(m_N^*)^2
=2(\beta_N-1)(1+o(1))$, while its numerator tends to one.
The susceptibility and fluctuation-scale asymptotics follow.
\end{proof}

\begin{remark}
The critical threshold is
\[
 p_N\asymp N^{-3/4},
\]
or equivalently $r_N=p_N^4N^3\asymp1$.  In particular, the effective inverse-temperature perturbation is negligible on the classical critical scale exactly when $r_N\to\infty$; Proposition~\ref{prop:annealed-critical-regimes}(1) makes this implication uniform, without requiring $p_N\to0$.
\end{remark}

\subsection{Annealed versus quenched critical crossover}

We first identify the annealed comparison density suggested by the
critical expansion at the beginning of this section.  In regime~(2) of
Proposition~\ref{prop:annealed-critical-regimes},
$r_N=p_N^4N^3\to c^{-2}$ and hence
$\sqrt N(\beta_N-1)\to c/24$.  On the critical magnetisation
scale $m=N^{-1/4}x$, the expansion of the effective Curie--Weiss
exponent derived above therefore has the limit
$cx^2/48-x^4/12$.  This motivates the annealed crossover target
\begin{equation}\label{eq:pc-def}
 p_c(x)=\frac1{Z_c}\exp\left\{\frac{c}{48}x^2-\frac{x^4}{12}\right\},
 \qquad c>0,
\end{equation}
where $Z_c$ is the normalising constant.  The same formula with $c=0$
gives the classical quartic critical target.  At this point the
expansion motivates the density; the quantitative convergence will be
proved below by Stein's method.

The quenched analysis identifies the same threshold $p_N^4N^3$, but the crossover law differs.  More precisely, \cite[Theorem~1.2]{KLS20} states that, if
\[
 p_N^4N^3\to\frac1{c^2},
\]
then the quenched $N^{-3/4}M_N$ limit has density
\begin{equation}\label{eq:KLS-crossover}
 f_c^{\rm que}(x)
 =\frac1{Z_c^{\rm que}}
 \exp\left\{-\frac{c}{24}x^2-\frac{x^4}{12}\right\}.
\end{equation}
By contrast, the annealed density \eqref{eq:pc-def} has quadratic
coefficient $+c/48$, whereas the quenched density
\eqref{eq:KLS-crossover} has coefficient $-c/24$.
Thus not only the coefficient but even the \emph{sign} of the quadratic term is different.  This is not a contradiction: annealing replaces the random graph by the deterministic effective Curie--Weiss parameter $\beta_N$, and for $p_N\to0$ the correction $\beta_N-1$ is positive.  In the quenched model, the disorder enters the limiting free-energy expansion through a different critical correction.

The third regime is equally different.  In \cite{KLS20}, it is shown that when $p_NN\to\infty$ but $p_N^4N^3\to0$, the quenched magnetisation satisfies
\[
 \frac{M_N}{p_NN^{3/2}}\Longrightarrow\mathcal N(0,12).
\]
For the annealed model, Proposition~\ref{prop:annealed-critical-regimes}(3) instead gives two supercritical wells whose centers, measured on the same $p_NN^{3/2}$ scale, diverge because
\[
 \frac{Nm_N^*}{p_NN^{3/2}}
 \sim\frac1{\sqrt8\,p_N^2N^{3/2}}
 =\frac1{\sqrt8\sqrt{p_N^4N^3}}
 \longrightarrow\infty.
\]
The susceptibility calculation identifies the corresponding local scale as $\sqrt{12}\,p_NN^{3/2}$.  We do not pursue a conditional Gaussian limit in this regime.  Hence the annealed and quenched models have the same quartic critical limiting law in the strongly dense regime $p_N^4N^3\to\infty$, but differ qualitatively on and below the dilution-critical line.

\subsection{Kolmogorov Stein theory for the annealed crossover law}
\label{subsec:crossover-Stein}

We derive bounds for the solutions of the indicator Stein equation associated with the crossover densities $p_c$, $c\ge0$, defined in \eqref{eq:pc-def}, to obtain a quantitative approximation in Kolmogorov distance.  Write $F_c$ for its distribution function and
$Q_c=1-F_c$ for the corresponding survival function.  Put
\[
 V_c(x)=\frac{x^4}{12}-\frac{c}{48}x^2,
 \qquad
 g_c(x)=V_c'(x)=\frac{x^3}{3}-\frac{c}{24}x.
\]
The density satisfies $p_c'/p_c=-g_c$, and hence the canonical Stein operator is
\[
 \mathcal A_cf(x)=f'(x)-g_c(x)f(x)
 =f'(x)+\left(\frac{c}{24}x-\frac{x^3}{3}\right)f(x).
\]
For $c=0$ this is the classical quartic Curie--Weiss operator; for $c>0$ the potential $V_c$ has two wells.  General representations of Stein solutions and bounds on the solutions and their derivatives can be found in \cite{LRS17,ERS20,ES22}; recent results for quartic and related densities appear in \cite{RapinSwan26,EicMills26}.  Here we derive explicit bounds for the indicator Stein equation associated with the nonconvex density $p_c$.

For related quartic targets in an edge--triangle model,
\cite[Lemma~3.1 and Section~8.1]{LSSX26} obtains bounds for
$1$-Lipschitz tests in Wasserstein distance.  In contrast, the
indicator tests required here are discontinuous, and their Stein
solutions and first derivatives are estimated directly below.

\begin{remark}\label{rem:EL10-positive-gamma}
In \cite[Theorem~1.4(1)]{EL10} we also claimed a Kolmogorov Berry--Esseen bound for the critical Curie--Weiss model with $\beta_N-1=\gamma N^{-1/2}$ and $\gamma>0$. However, the proof does not establish this case: the bounds for Stein solutions invoked there were verified for unimodal pure-power targets, not for the double-well density proportional to $\exp\{-x^4/12+\gamma x^2/2\}$. The present analysis fills this gap by establishing the required bounds for the indicator Stein solutions and their first derivatives.
\end{remark}

For the indicator test $H_z(x)=\1_{\{x\le z\}}$ introduced
in Section~\ref{sec:introduction}, the bounded Stein solution is
\begin{equation}\label{eq:fz-crossover}
 f_{c,z}(x)=
 \begin{cases}
 \displaystyle \frac{F_c(x)Q_c(z)}{p_c(x)},&x\le z,\\[2mm]
 \displaystyle \frac{F_c(z)Q_c(x)}{p_c(x)},&x>z.
 \end{cases}
\end{equation}
We shall write
\[
 \ell_c(x)=\frac{F_c(x)}{p_c(x)},
 \qquad
 r_c(x)=\frac{Q_c(x)}{p_c(x)}.
\]
Then
\begin{equation}\label{eq:ell-r-ode}
 \ell_c'(x)=1+g_c(x)\ell_c(x),
 \qquad
 r_c'(x)=g_c(x)r_c(x)-1.
\end{equation}
These two one-dimensional functions contain the entire first-order Kolmogorov Stein analysis.

For fixed $x$, maximizing \eqref{eq:fz-crossover} over $z$ gives
\[
 \sup_{z\in\R}|f_{c,z}(x)|
 =K_{0,c}(x):=\frac{F_c(x)Q_c(x)}{p_c(x)}.
\]
The following proposition determines the global maximum.

\begin{proposition}[Optimal bound on indicator Stein solutions]\label{prop:CK0-crossover}
For every $c\ge0$,
\[
 K_0(c):=\sup_{z\in\R}\|f_{c,z}\|_\infty
 =K_{0,c}(0)=\frac{1}{4p_c(0)}=\frac{Z_c}{4}.
\]
In particular, this bound is finite for every fixed $c$.
\end{proposition}

\begin{proof}
By symmetry, $K_{0,c}$ is even.  Put
$q_c(x)=\exp\{cx^2/48-x^4/12\}$, so that $p_c=q_c/Z_c$, and
$d_c=\sqrt c/2$.  For $0\le x\le d_c$ we have
$cx^2/48-x^4/12=x^2(c-4x^2)/48\ge0$; hence $q_c(x)\ge1$ and
$K_{0,c}(x)\le Z_c/(4q_c(x))\le Z_c/4$, since $F_cQ_c\le1/4$.

For $x\ge d_c$, set $y(x)=2F_c(x)-1$ and
$R(x)=y(x)^2+q_c(x)-1$.  Because $F_cQ_c=(1-y^2)/4$,
$K_{0,c}(x)\le Z_c/4$ is equivalent to $R(x)\ge0$.
Using $y'=2q_c/Z_c$ and $q_c'=-g_cq_c$, we obtain
\[
 R'(x)=q_c(x)J(x),\qquad
 J(x)=\frac{4y(x)}{Z_c}-g_c(x),\qquad
 J'(x)=\frac{8q_c(x)}{Z_c^2}-g_c'(x).
\]
Moreover,
$J''(x)=-8g_c(x)q_c(x)/Z_c^2-2x<0$ for $x>d_c$,
while $R(d_c)=y(d_c)^2\ge0$, $R(x)\to0$ and
$J(x)\to-\infty$ as $x\to\infty$.

If $c=0$, then $d_c=0$, $J(0)=0$ and $J'(0)=8/Z_0^2>0$.
Strict concavity on $(0,\infty)$ shows that $J$ is first positive
and then negative, with exactly one positive zero.  Thus $R$
first increases from $R(0)=0$ and then decreases to $0$.

Now let $c>0$.  If $J(d_c)\ge0$, strict concavity likewise shows
that $J$ can change sign only once, from positive to negative;
therefore $R$ first increases and then decreases to $0$.
If $J(d_c)<0$, write $A_c=\int_0^{d_c}q_c(t)\,dt\ge d_c$.
Since $y(d_c)=2A_c/Z_c$ and $g_c(d_c)=c^{3/2}/48$,
the inequality $J(d_c)<0$ implies
$Z_c^2>384A_c/c^{3/2}\ge192/c$.
As $q_c(d_c)=1$ and $g_c'(d_c)=5c/24$, it follows that
\[
 J'(d_c)=\frac8{Z_c^2}-\frac{5c}{24}< -\frac c6<0.
\]
By strict concavity, $J'<0$ and $J<0$ on $[d_c,\infty)$;
therefore $R$ decreases to $0$ and is nonnegative there.
Together with the first interval, this proves the proposition.
\end{proof}

The next proposition gives the uniform bound on the first derivatives of the indicator Stein solutions needed for the quantitative Kolmogorov approximation.  Derivatives are understood almost everywhere, with the essential supremum norm; the same bounds hold for the one-sided derivatives at $x=z$.

\begin{proposition}[Bounds for indicator Stein solutions and their derivatives]\label{prop:crossover-simple-factors}
For every $c\ge0$ and $z\in\R$, the solution $f_{c,z}$ in \eqref{eq:fz-crossover} satisfies
\[
 \|f_{c,z}\|_\infty\le K_0(c):=Z_c/4
\]
and
\begin{equation}\label{eq:simple-K1}
 \|f'_{c,z}\|_\infty\le K_1(c):=
 1+\frac{Z_c}{6}\left(\frac c{24}\right)^{3/2}.
\end{equation}
The first bound is optimal by Proposition~\ref{prop:CK0-crossover}.
\end{proposition}

\begin{proof}
The bound for $f_{c,z}$ is Proposition~\ref{prop:CK0-crossover}.
The Stein equation gives $|f'_{c,z}(x)|\le1+|g_c(x)|\,|f_{c,z}(x)|$.
Set $a_c=\sqrt{c/8}$.  On $0\le x\le a_c$ we have $g_c(x)\le0$ and
$|g_c(x)|=cx/24-x^3/3$, whose maximum, attained at $x=\sqrt{c/24}$, is $\frac23(c/24)^{3/2}$.  Hence
\[
 |f'_{c,z}(x)|\le1+\frac{Z_c}{6}\left(\frac c{24}\right)^{3/2}
 \qquad(0\le x\le a_c).
\]
For $x>a_c$, both $g_c(x)$ and $g_c'(x)$ are positive.  Thus $g_c(t)\ge g_c(x)$ for $t\ge x$ and, using $p_c'=-g_cp_c$,
\[
 Q_c(x)=\int_x^\infty p_c(t)\,dt
 \le\frac1{g_c(x)}\int_x^\infty g_c(t)p_c(t)\,dt
 =\frac{p_c(x)}{g_c(x)}.
\]
Consequently $0\le-r_c'(x)=1-g_c(x)r_c(x)\le1$.  The identities in \eqref{eq:ell-r-ode} also give
$Q_c(x)\ell_c'(x)+F_c(x)(-r_c'(x))=1$, and both summands are nonnegative.  Hence $0\le Q_c(x)\ell_c'(x)\le1$.  In \eqref{eq:fz-crossover}, if $x\le z$ then $Q_c(z)\le Q_c(x)$, while if $x>z$ then $F_c(z)\le F_c(x)$.  In either case $|f'_{c,z}(x)|\le1$ for $x>a_c$.  Finally, symmetry yields $f_{c,z}(-x)=f_{c,-z}(x)$, so the same estimates hold for $x<0$.  This proves \eqref{eq:simple-K1}.
\end{proof}

\begin{remark}\label{rem:optional-sharpening}
The derivative bound \eqref{eq:simple-K1} is deliberately not optimized.  The optimal uniform derivative bound can be expressed by a one-dimensional variational formula involving $F_c/p_c$ and $(1-F_c)/p_c$, and its sharp value exhibits a nontrivial transition as $c$ varies.  This refinement is analytically interesting but plays no role in the $N^{-1/2}$ Kolmogorov argument; it is recorded separately in Appendix~\ref{app:optimal-K1}.
\end{remark}

\begin{remark}
The normalizing constant can also be removed from the preceding bounds.  The elementary inequality
\[
 \frac{c}{48}x^2-\frac{x^4}{12}
 \le \frac{c^2}{384}-\frac{x^4}{24}
\]
gives
\[
 Z_c\le \frac{24^{1/4}}2\Gamma\!\left(\frac14\right)e^{c^2/384}.
\]
Combining this estimate with Proposition~\ref{prop:crossover-simple-factors} gives explicit constants depending only on $c$.
\end{remark}

\subsection{A Kolmogorov plug-in theorem}
\label{subsec:crossover-plugin}

We first establish the general bound that will determine which quantities need to be controlled in the annealed Curie--Weiss model. It uses the solutions of the indicator Stein equations from Proposition~\ref{prop:crossover-simple-factors}.

\begin{theorem}[Kolmogorov plug-in bound]\label{prop:crossover-plugin}
Let $(W,W')$ be exchangeable, $D=W-W'$, let $\lambda>0$ and $\delta>0$, and suppose
\[
 |D|\le\delta,
 \qquad
 \E(D\mid W)=\lambda\{g_c(W)+R(W)\},
\]
with
\[
 A:=\E\left|1-\frac{\E(D^2\mid W)}{2\lambda}\right|<\infty,
 \qquad
 D^2\le4\lambda\quad\text{a.s.}
\]
Assume $\E W^2+\E|W|^3<\infty$.  Then
\begin{align}\label{eq:crossover-plugin}
 d_K(W,Y_c)
 \le{}& K_0(c)\E|R(W)|
 +(K_1(c)+2)A\\
 &+C_c\delta\{1+\E W^2+\E|W|^3\},\nonumber
\end{align}
where $Y_c$ has density \eqref{eq:pc-def}, and the explicit bounds $K_0(c),K_1(c)$ for the Stein solution and its derivative are given in Proposition~\ref{prop:crossover-simple-factors}.  The constant $C_c<\infty$ depends only on the fixed target parameter $c$.
\end{theorem}

\begin{proof}
Write $g=g_c$, $F=F_c$, $p=p_c$, and let $f=f_{c,z}$ be the bounded
solution of the indicator Stein equation
\begin{equation}\label{eq:crossover-Stein-indicator-proof}
 f'(w)-g(w)f(w)=\1_{\{w\le z\}}-F(z).
\end{equation}
In standard exchangeable-pair arguments for Kolmogorov bounds,
monotonicity of the Stein drift is often used to control the
discontinuity of the indicator test function.  For $c>0$, this
route is unavailable because $g_c$ is nonmonotone.  We therefore
use the positivity of the exchangeable-pair jump kernel to obtain
separate one-sided estimates; the argument also applies when $c=0$.

If $\delta>1$, then \eqref{eq:crossover-plugin} is trivial after increasing
$C_c$, since $d_K\le1$ and the last factor on the right hand side is at
least $C_c\delta$.  We may therefore assume throughout that $0<\delta\le1$.

Set
\[
 K(t)=\E\left[D\left\{\1_{\{-D\le t\le0\}}-
                 \1_{\{0<t\le-D\}}\right\}\,\middle|\,W\right].
\]
For every realization of $D$ the quantity inside the conditional
expectation is nonnegative; hence $K(t)\ge0$.  Moreover $K(t)=0$ for
$|t|>\delta$ and
\begin{equation}\label{eq:crossover-kernel-mass}
 \int_{-\delta}^{\delta}K(t)\,dt=\E(D^2\mid W).
\end{equation}
Exchangeability gives, for every absolutely continuous $f$ for which the
expectations exist,
\begin{align}
 2\E[D f(W)]
 &=\E\{D(f(W)-f(W'))\}\notag\\
 &=\E\int_{-\delta}^{\delta} f'(W+t)K(t)\,dt .
 \label{eq:crossover-kernel-identity-full}
\end{align}
Using the regression assumption, the left hand side is
$2\lambda\E[f(W)\{g(W)+R(W)\}]$.

We first prove the upper one-sided bound.  Since $K\ge0$ and $|t|\le\delta$
on its support,
\[
 \1_{\{W+t\le z\}}\ge \1_{\{W\le z-\delta\}}.
\]
Insert \eqref{eq:crossover-Stein-indicator-proof} into
\eqref{eq:crossover-kernel-identity-full}.  Then
\begin{align}
 &\E\!\left[\1_{\{W\le z-\delta\}}D^2\right]
     -F(z)\E D^2                                           \notag\\
 &\quad\le
 2\lambda\E[f(W)g(W)]+2\lambda\E[f(W)R(W)]
 -\E\int_{-\delta}^{\delta}g(W+t)f(W+t)K(t)\,dt            \notag\\
 &\quad=
 2\lambda\E\!\left[f(W)g(W)
 \left\{1-\frac{\E(D^2\mid W)}{2\lambda}\right\}\right]
 +2\lambda\E[f(W)R(W)]                                    \\
 &\qquad\quad+
 \E\int_{-\delta}^{\delta}
 \{g(W)f(W)-g(W+t)f(W+t)\}K(t)\,dt                         \notag\\
 &\quad=:J_1+J_2+J_3 .
 \label{eq:crossover-J123}
\end{align}
The Stein equation and Proposition~\ref{prop:crossover-simple-factors}
give
\[
 \|gf\|_\infty
 \le \|f'\|_\infty+1
 \le K_1(c)+1.
\]
Consequently,
\begin{equation}\label{eq:J1J2-crossover}
 |J_1|\le2\lambda\{K_1(c)+1\}A,
 \qquad
 |J_2|\le2\lambda K_0(c)\E|R(W)|.
\end{equation}

It remains to control $J_3$.  Put $q(w)=g(w)f(w)$.  For $|t|\le\delta$,
\begin{align}
 |q(w+t)-q(w)|
 &\le |g(w)|\,|f(w+t)-f(w)|
      +|f(w+t)|\,|g(w+t)-g(w)|\notag\\
 &\le |t|\left\{K_1(c)|g(w)|
      +K_0(c)\sup_{|s-w|\le\delta}|g'(s)|\right\}.
 \notag
\end{align}
For the polynomial
\[
 g(w)=\frac{w^3}{3}-\frac c{24}w,
 \qquad g'(w)=w^2-\frac c{24},
\]
and $\delta\le1$, there is a finite constant $L_c$ such that
\[
 |g(w)|\le L_c(1+|w|^3),
 \qquad
 \sup_{|s-w|\le\delta}|g'(s)|\le L_c(1+w^2).
\]
Thus, after increasing $L_c$ if necessary,
\[
 |q(w+t)-q(w)|
 \le L_c |t|\{1+w^2+|w|^3\}.
\]
Using $K\ge0$, $|t|\le\delta$, \eqref{eq:crossover-kernel-mass}, and
$D^2\le4\lambda$ almost surely, we obtain
\begin{align}
 |J_3|
 &\le L_c\delta\,
 \E\left[\{1+W^2+|W|^3\}\E(D^2\mid W)\right]\notag\\
 &\le4\lambda L_c\delta
 \{1+\E W^2+\E|W|^3\}.
 \label{eq:J3-crossover}
\end{align}

We now convert the left hand side of \eqref{eq:crossover-J123} into a
Kolmogorov difference.  Set
\[
 V(W)=\frac{\E(D^2\mid W)}{2\lambda}.
\]
Then, using $|\1_{\{W\le z-\delta\}}-F(z)|\le1$,
\begin{align}
 &\E\!\left[\1_{\{W\le z-\delta\}}D^2\right]-F(z)\E D^2
 \notag\\
 &=2\lambda\{\P(W\le z-\delta)-F(z)\}
 -2\lambda\E\!\left[
   \{\1_{\{W\le z-\delta\}}-F(z)\}\{1-V(W)\}\right]
 \notag\\
 &\ge2\lambda\{\P(W\le z-\delta)-F(z-\delta)\}
       -2\lambda A
       -2\lambda\{F(z)-F(z-\delta)\}                         \notag\\
 &\ge2\lambda\{\P(W\le z-\delta)-F(z-\delta)\}
       -2\lambda A-2\lambda\|p_c\|_\infty\delta .
 \label{eq:left-lower-crossover}
\end{align}
Combining \eqref{eq:crossover-J123}--\eqref{eq:left-lower-crossover} and
dividing by $2\lambda$ gives
\begin{align}
 \P(W\le z-\delta)-F(z-\delta)
 \le{}&\{K_1(c)+2\}A+K_0(c)\E|R(W)|\notag\\
 &+C_c\delta\{1+\E W^2+\E|W|^3\},
 \label{eq:upper-one-sided-crossover}
\end{align}
where $\|p_c\|_\infty\delta$ has been absorbed into the last term.
Since $z-\delta$ ranges over all of $\R$, this proves the required upper
one-sided Kolmogorov bound.

For the lower one-sided bound use instead
\[
 \1_{\{W+t\le z\}}\le \1_{\{W\le z+\delta\}}.
\]
The same substitution in \eqref{eq:crossover-kernel-identity-full} now gives
\begin{equation}\label{eq:crossover-lower-kernel-raw}
 \E\!\left[\1_{\{W\le z+\delta\}}D^2\right]-F(z)\E D^2
 \ge J_1+J_2+J_3,
\end{equation}
with the same $J_1,J_2,J_3$ as in \eqref{eq:crossover-J123}; hence the
absolute-value bounds \eqref{eq:J1J2-crossover} and
\eqref{eq:J3-crossover} imply that the left hand side of
\eqref{eq:crossover-lower-kernel-raw} is bounded below by
\begin{equation}\label{eq:crossover-lower-kernel-bound}
 -2\lambda\Big[\{K_1(c)+1\}A+K_0(c)\E|R(W)|
 +C_c\delta\{1+\E W^2+\E|W|^3\}\Big].
\end{equation}
On the other hand, with $x=z+\delta$ and the same $V(W)$ as above,
\begin{align}
 &\E\!\left[\1_{\{W\le x\}}D^2\right]-F(z)\E D^2
 \notag\\
 &=2\lambda\{\P(W\le x)-F(z)\}
 -2\lambda\E\!\left[
   \{\1_{\{W\le x\}}-F(z)\}\{1-V(W)\}\right]
 \notag\\
 &\le2\lambda\{\P(W\le x)-F(x)\}
       +2\lambda\{F(x)-F(z)\}+2\lambda A\notag\\
 &\le2\lambda\{\P(W\le x)-F(x)\}
       +2\lambda\|p_c\|_\infty\delta+2\lambda A.
 \label{eq:left-upper-crossover}
\end{align}
Combining \eqref{eq:crossover-lower-kernel-bound} with
\eqref{eq:left-upper-crossover} and dividing by $2\lambda$ yields
\begin{align}
 F(z+\delta)-\P(W\le z+\delta)
 \le{}&\{K_1(c)+2\}A+K_0(c)\E|R(W)|\notag\\
 &+C_c\delta\{1+\E W^2+\E|W|^3\}.
 \label{eq:lower-one-sided-crossover}
\end{align}
Again $z+\delta$ ranges over $\R$.  Taking the larger of
\eqref{eq:upper-one-sided-crossover} and
\eqref{eq:lower-one-sided-crossover} proves
\eqref{eq:crossover-plugin}.

We emphasize that the proof used only exchangeability, the bounded jump,
the regression and conditional-variance assumptions, the explicit bounds
on $f_{c,z}$ and $f_{c,z}'$, and the polynomial growth of $g_c,g_c'$.  No
monotonicity or sign condition on $g_c$ enters the argument.
\end{proof}

\subsection{Verification for the annealed Curie--Weiss model}
\label{subsec:critical-exchangeable}

We now verify the regression, conditional-variance and moment assumptions of Theorem~\ref{prop:crossover-plugin} for the annealed Curie--Weiss model.

Under the exact annealed reduction \eqref{eq:intro-exact-CW}, at $h=0$ the spin law is exactly $CW(\beta_N,0)$.  Choose $I$ uniformly from $\{1,\ldots,N\}$ and resample $\sigma_I'$ from its conditional Gibbs law given all remaining spins.  Set
\begin{equation}\label{eq:critical-W-pair}
 W_N=\frac{M_N}{N^{3/4}},
 \qquad
 W_N'=\frac{M_N-\sigma_I+\sigma_I'}{N^{3/4}},
 \qquad
 D_N=W_N-W_N'.
\end{equation}
The one-step heat-bath dynamics satisfies detailed balance with respect to $CW(\beta_N,0)$ and is therefore reversible; hence $(W_N,W_N')$ is exchangeable, and
$|D_N|\le2N^{-3/4}$.
Put $\theta_N=\sqrt N(\beta_N-1)$ and $\lambda_N=N^{-3/2}$.

\begin{lemma}[Critical drift regression]\label{lem:critical-drift}
Fix $c\ge0$ and $K<\infty$. Suppose $|\theta_N|\le K$, where
$\beta_N=1+\theta_NN^{-1/2}$, and let $W_N$ be as above. Then, for all sufficiently large $N$,
\begin{equation}\label{eq:critical-regression}
 \E(D_N\mid W_N)=\lambda_N\{g_c(W_N)+R_{N,c}(W_N)\},
 \qquad \lambda_N=N^{-3/2},
\end{equation}
where, uniformly for all admissible values $|w|\le N^{1/4}$,
\begin{align}\notag
 R_{N,c}(w)
 ={}&-\left(\theta_N-\frac c{24}\right)w
 +\frac{\beta_N^3-1}{3}w^3+\rho_{N,c}(w),\\
 |\rho_{N,c}(w)|
 \le{}& C_{c,K}N^{-1/2}\bigl(1+|w|+|w|^5\bigr).
 \label{eq:Rcritical-rho}
\end{align}
Consequently, if $\sup_N\E|W_N|^5<\infty$, then
\begin{equation}\label{eq:Rcritical-L1}
 \E|R_{N,c}(W_N)|
 \le C_{c,K}\left\{N^{-1/2}+\left|\theta_N-\frac c{24}\right|\right\}.
\end{equation}
\end{lemma}

\begin{proof}
Write
\[
 m=\frac{M_N}{N}=N^{-1/4}W_N,
 \qquad T(x)=\tanh x,
 \qquad t=\beta_Nm,
 \qquad \varepsilon_N=\frac{\beta_N}{N}.
\]
Since $|\theta_N|\le K$, we may assume $1/2\le \beta_N\le2$ for all sufficiently large $N$.
Averaging the heat-bath conditional expectation over the uniformly chosen site gives the exact identity
\[
 \E(D_N\mid W_N)=N^{-3/4}\{m-A_N(m)\},
\]
where
\[
 A_N(m)=\frac{1+m}{2}T(t-\varepsilon_N)
       +\frac{1-m}{2}T(t+\varepsilon_N).
\]
We first control the finite-site shift without destroying the symmetry at $m=0$.
Taylor's formula to second order gives, uniformly for $|t|\le2$,
\[
 T(t\pm\varepsilon_N)
 =T(t)\pm\varepsilon_NT'(t)
   +\frac{\varepsilon_N^2}{2}T''(t)+O_K(\varepsilon_N^3).
\]
Since $T''(0)=0$ and $T'''$ is bounded on $[-3,3]$,
\[
 |T''(t)|\le C|t|\le C_K|m|.
\]
Substitution into the preceding convex combination yields
\begin{equation}\label{eq:A-critical-shift-sharp}
 A_N(m)
 =T(\beta_Nm)-\frac{\beta_Nm}{N}T'(\beta_Nm)+r_{A,N}(m),
 \qquad
 |r_{A,N}(m)|\le C_K\left(\frac{|m|}{N^2}+\frac1{N^3}\right).
\end{equation}

Taylor's theorem at the origin, with a seventh-order remainder uniform for $|\beta_Nm|\le2$, gives
\begin{equation}\label{eq:tanh-critical-seventh}
 T(\beta_Nm)
 =\beta_Nm-\frac{\beta_N^3m^3}{3}+\frac{2\beta_N^5m^5}{15}+r_{7,N}(m),
 \qquad |r_{7,N}(m)|\le C_K|m|^7.
\end{equation}
Using \eqref{eq:A-critical-shift-sharp}, \eqref{eq:tanh-critical-seventh} and $|T'|\le1$ in the preceding exact drift identity, we obtain
\begin{align*}
 \E(D_N\mid W_N)
 ={}&N^{-3/4}\left\{(1-\beta_N)m+\frac{\beta_N^3}{3}m^3\right\}
       +\mathcal R_N(w)\\
 ={}&N^{-3/2}\left\{-\theta_Nw+\frac{\beta_N^3}{3}w^3\right\}
       +\mathcal R_N(w),
\end{align*}
where, because $|w|\le N^{1/4}$,
\begin{align*}
 |\mathcal R_N(w)|
 &\le C_K\left(N^{-2}|w|^5+N^{-5/2}|w|^7+N^{-2}|w|
 +N^{-3}|w|+N^{-15/4}\right)\\
 &\le C_KN^{-2}\bigl(1+|w|+|w|^5\bigr).
\end{align*}
Here $N^{-5/2}|w|^7\le N^{-2}|w|^5$ because $|w|\le N^{1/4}$.
Adding and subtracting $cw/24$ and $w^3/3$, and dividing the remainder by $\lambda_N=N^{-3/2}$, proves \eqref{eq:critical-regression}--\eqref{eq:Rcritical-rho}.
Finally,
\[
 |\beta_N^3-1|=|\beta_N-1|(\beta_N^2+\beta_N+1)\le C_KN^{-1/2}.
\]
Taking expectations and using the fifth-moment bound gives \eqref{eq:Rcritical-L1}.
\end{proof}

\begin{lemma}[Conditional variance]\label{lem:critical-variance}
Fix $K<\infty$ and suppose $|\theta_N|\le K$. Then, for all sufficiently large $N$,
\[
 \frac{\E(D_N^2\mid W_N)}{2\lambda_N}=1+\Delta_N(W_N),
\]
where, uniformly for all admissible $|w|\le N^{1/4}$,
\begin{equation}\label{eq:Scritical-pointwise}
 \Delta_N(w)=-\beta_NN^{-1/2}w^2+\eta_N(w),
 \qquad |\eta_N(w)|\le C_KN^{-1}(1+w^4).
\end{equation}
Hence, if $\sup_N\E W_N^4<\infty$,
\begin{equation}\label{eq:Scritical-L1}
 \E|\Delta_N(W_N)|\le C_KN^{-1/2}.
\end{equation}
Moreover
\begin{equation}\label{eq:critical-D2-bound}
 D_N^2\le4N^{-3/2}=4\lambda_N\qquad\text{a.s.}
\end{equation}
\end{lemma}

\begin{proof}
Again write $m=N^{-1/4}W_N$, $T(x)=\tanh x$, $t=\beta_Nm$ and $\varepsilon_N=\beta_N/N$.
Since $(\sigma_i-\sigma_i')^2=2(1-\sigma_i\sigma_i')$, averaging over the resampled site gives
\[
 \E(D_N^2\mid W_N)=2N^{-3/2}\{1-B_N(m)\},
\]
with
\[
 B_N(m)=\frac{1+m}{2}T(t-\varepsilon_N)
       -\frac{1-m}{2}T(t+\varepsilon_N).
\]
Taylor expansion in $\varepsilon_N$ gives
\[
 B_N(m)=mT(\beta_Nm)-\frac{\beta_N}{N}T'(\beta_Nm)+r_{B,N}(m),
 \qquad |r_{B,N}(m)|\le C_KN^{-2}.
\]
Indeed, the second-order contribution is $m\varepsilon_N^2T''(t)/2=O_K(N^{-2})$, and the third-order remainder is $O_K(N^{-3})$.
On $|\beta_Nm|\le2$,
\[
 T(\beta_Nm)=\beta_Nm-\frac{\beta_N^3m^3}{3}+O_K(|m|^5),
 \qquad T'(\beta_Nm)=1+O_K(m^2).
\]
Therefore
\begin{align*}
 B_N(m)
 ={}&\beta_Nm^2-\frac{\beta_N^3}{3}m^4-\frac{\beta_N}{N}
 +O_K\left(m^6+\frac{m^2}{N}+\frac1{N^2}\right)\\
 ={}&\beta_NN^{-1/2}w^2+\widetilde\eta_N(w),
\end{align*}
where
\begin{align*}
 |\widetilde\eta_N(w)|
 &\le C_K\left(N^{-1}w^4+N^{-3/2}|w|^6+N^{-1}
 +N^{-3/2}w^2+N^{-2}\right)\\
 &\le C_KN^{-1}(1+w^4).
\end{align*}
Here we used $|w|\le N^{1/4}$, so that
$N^{-3/2}|w|^6\le N^{-1}w^4$ and $N^{-3/2}w^2\le N^{-1}$.
Since the preceding identity and $2\lambda_N=2N^{-3/2}$ give
\[
 \frac{\E(D_N^2\mid W_N)}{2\lambda_N}=1-B_N(m),
\]
we take $\Delta_N(w)=-B_N(N^{-1/4}w)$, proving \eqref{eq:Scritical-pointwise}. The fourth-moment bound gives \eqref{eq:Scritical-L1}, and the bounded-jump estimate above gives \eqref{eq:critical-D2-bound}.
\end{proof}

\begin{lemma}[Uniform critical moments]\label{lem:critical-moments}
Fix $K<\infty$ and suppose $\beta_N=1+\theta_NN^{-1/2}$, with
$|\theta_N|\le K$. Under $CW(\beta_N,0)$, the magnetization
$W_N=N^{-3/4}M_N$ satisfies
\begin{equation}\label{eq:critical-moments}
 \sup_{N}\E|W_N|^5<\infty.
\end{equation}
The bound depends only on $K$.
\end{lemma}

\begin{proof}
We use the exchangeable pair defined in \eqref{eq:critical-W-pair}.
Write $m=N^{-1/4}w$, $t=\beta_Nm$ and $\varepsilon_N=\beta_N/N$.
The exact drift identity from Lemma~\ref{lem:critical-drift} holds with
$A_N(m)$ as defined there. Since $\beta_N\in[1/2,2]$ for large $N$,
the symmetric average of $\tanh(t\pm\varepsilon_N)$ differs from
$\tanh(t)$ by at most $C\varepsilon_N^2|t|$, whereas the difference
of the two terms is at most $2\varepsilon_N$. Consequently
$|A_N(m)-\tanh(\beta_Nm)|\le C|m|/N$ uniformly for $|m|\le1$.
Here we used $|\tanh''(v)|\le C|v|$ on a fixed bounded interval.

There is an $a>0$, independent of $N$, such that
$x\{x-\tanh x\}\ge a x^4$ for $|x|\le2$; this follows from
$\lim_{x\to0}(x-\tanh x)/x^3=1/3$ and positivity away from zero.
With $t=\beta_Nm$, the identity
\[
 m\{m-\tanh(\beta_Nm)\}
 =(1-\beta_N)m^2+\beta_N^{-1}t\{t-\tanh t\}
 \ge \frac a8 m^4-KN^{-1/2}m^2
\]
holds uniformly for $|m|\le1$ and all sufficiently large $N$:
indeed, $|t|\le2$, $\beta_N\ge1/2$, and $\beta_N^3\ge1/8$.
The finite-site bound above implies
$|m\{A_N(m)-\tanh(\beta_Nm)\}|\le Cm^2/N$.
Since $w=N^{1/4}m$, $\lambda_N=N^{-3/2}$, and the exact
conditional drift equals $N^{-3/4}\{m-A_N(m)\}$, we obtain,
uniformly over all admissible $w$,
\begin{equation}\label{eq:critical-drift-coercive}
 w\E(D_N\mid W_N=w)
 \ge \lambda_N\{a_0w^4-C_Kw^2\},\qquad a_0=a/8>0.
\end{equation}

By exchangeability and $D_N=W_N-W_N'$,
$\E[D_NW_N]=\frac12\E D_N^2\le2\lambda_N$.
Combining this identity with \eqref{eq:critical-drift-coercive} gives
$a_0\E W_N^4\le 2+C_K\E W_N^2$.
Since $w^2\le\epsilon w^4+C_\epsilon$, choosing $\epsilon$ small
shows that $\sup_N\E W_N^4<\infty$.

To obtain the sixth moment, use the same exchangeable-pair identity
with the function $w\mapsto w^3$:
\[
 \E[D_NW_N^3]
 =\frac12\E\bigl[D_N^2\{3W_N^2-3W_ND_N+D_N^2\}\bigr].
\]
The bounds $|D_N|\le2N^{-3/4}\le2$ and
$D_N^2\le4\lambda_N$ imply
$|\E[D_NW_N^3]|\le C\lambda_N(1+\E W_N^2)$.
Multiplying \eqref{eq:critical-drift-coercive} by $w^2$ and taking
expectations now yields
$a_0\E W_N^6\le C(1+\E W_N^2)+C_K\E W_N^4$.
Thus $\sup_N\E W_N^6<\infty$, and H\"older's inequality gives
\eqref{eq:critical-moments}. The finitely many smaller $N$ are harmless.
\end{proof}

\begin{theorem}[Kolmogorov rate in the annealed crossover window]\label{thm:crossover-K}
Fix $c\ge0$ and assume $h=0$, nominal $\beta=1$, $p_NN\to\infty$, and
\[
 \theta_N:=\sqrt N(\beta_N-1)=\frac c{24}+o(1).
\]
Let $W_N$ be defined by \eqref{eq:critical-W-pair}, and let $Y_c$ have density \eqref{eq:pc-def}.  Then, for all sufficiently large $N$,
\begin{equation}\label{eq:crossover-K-rate}
 d_K(W_N,Y_c)
 \le C(c)\left\{N^{-1/2}+\left|\theta_N-\frac c{24}\right|\right\}.
\end{equation}
In particular, if
\begin{equation}\label{eq:theta-rate}
 \left|\sqrt N(\beta_N-1)-\frac c{24}\right|=O(N^{-1/2}),
\end{equation}
then
\begin{equation}\label{eq:crossover-opt-rate}
 d_K(W_N,Y_c)=O_c(N^{-1/2}).
\end{equation}
The result is valid for every fixed $c\ge0$.
\end{theorem}

\begin{proof}
We verify the hypotheses of Theorem~\ref{prop:crossover-plugin} one by one.
Since
\[
 \theta_N=\frac c{24}+o(1),
\]
there is, for example,
\[
 K_c:=1+\frac c{24}
\]
such that $|\theta_N|\le K_c$ for all sufficiently large $N$.  Thus all constants in Lemmas~\ref{lem:critical-drift}--\ref{lem:critical-moments} may be chosen depending only on the fixed target parameter $c$.

For the heat-bath pair \eqref{eq:critical-W-pair}, exchangeability was noted above and
\[
 |D_N|\le \delta_N:=2N^{-3/4}.
\]
Moreover, Lemma~\ref{lem:critical-drift} gives the regression identity in exactly the normalization required by Theorem~\ref{prop:crossover-plugin},
\[
 \E(D_N\mid W_N)
 =\lambda_N\{g_c(W_N)+R_{N,c}(W_N)\},
 \qquad \lambda_N=N^{-3/2},
\]
and, using the fifth-moment bound from Lemma~\ref{lem:critical-moments},
\begin{equation}\label{eq:crossover-assembly-R}
 \E|R_{N,c}(W_N)|
 \le C_c\left\{N^{-1/2}+\left|\theta_N-\frac c{24}\right|\right\}.
\end{equation}

Next, Lemma~\ref{lem:critical-variance} gives
\[
 A_N:=\E\left|1-\frac{\E(D_N^2\mid W_N)}{2\lambda_N}\right|
 \le C_cN^{-1/2},
\]
and also
\[
 D_N^2\le4\lambda_N\qquad\text{a.s.}
\]
Finally, Lemma~\ref{lem:critical-moments} and H\"older's inequality yield
\begin{equation}\label{eq:crossover-assembly-moments}
 \E W_N^2+\E|W_N|^3\le C_c.
\end{equation}
Therefore Theorem~\ref{prop:crossover-plugin} applies and gives
\begin{align*}
 d_K(W_N,Y_c)
 \le{}&K_0(c)\E|R_{N,c}(W_N)|
       +(K_1(c)+2)A_N\\
 &+C_c\delta_N\{1+\E W_N^2+\E|W_N|^3\}.
\end{align*}
Using \eqref{eq:crossover-assembly-R}--\eqref{eq:crossover-assembly-moments} and $\delta_N=2N^{-3/4}$, we obtain
\begin{align*}
 d_K(W_N,Y_c)
 &\le C_c\left\{N^{-1/2}+\left|\theta_N-\frac c{24}\right|\right\}
     +C_cN^{-3/4}\\
 &\le C_c'\left\{N^{-1/2}+\left|\theta_N-\frac c{24}\right|\right\},
\end{align*}
which is \eqref{eq:crossover-K-rate}.  If in addition \eqref{eq:theta-rate} holds, the right hand side is $O_c(N^{-1/2})$, proving \eqref{eq:crossover-opt-rate}.
\end{proof}

\begin{remark}\label{rem:dilution-rate}
If $p_N\to0$, Lemma~\ref{lem:critical-betaeff} gives
\[
 \theta_N
 =\frac{1+o(1)}{24\sqrt{p_N^4N^3}}.
\]
Thus $p_N^4N^3\to c^{-2}$ identifies the target $Y_c$, but a pure $O(N^{-1/2})$ Kolmogorov rate requires a corresponding quantitative rate of convergence of $\theta_N$ to $c/24$.  Without such a rate, \eqref{eq:crossover-K-rate} is the natural quantitative formulation delivered by the theorem: the probabilistic Stein error is $O_c(N^{-1/2})$, while the deterministic parameter mismatch is displayed separately.
\end{remark}

\begin{corollary}[Annealed critical distributional limits]\label{cor:annealed-critical-limits}
Assume $h=0$, nominal $\beta=1$ and $p_NN\to\infty$.
\begin{enumerate}
\item If $r_N=p_N^4N^3\to\infty$, then
\begin{equation}\label{eq:crit-quartic-limit}
 \frac{M_N}{N^{3/4}}\Longrightarrow Y_0,
\end{equation}
where $Y_0$ has density
\[
 f_0(x)=\frac1{Z_0}\exp\left(-\frac{x^4}{12}\right).
\]
More quantitatively, by \eqref{eq:theta-full-expansion},
\begin{equation}\label{eq:strong-dense-critical-rate}
 d_K\!\left(\frac{M_N}{N^{3/4}},Y_0\right)
 \le C\left\{N^{-1/2}+r_N^{-1/2}+(r_N\sqrt N)^{-1}\right\}
\end{equation}
for all sufficiently large $N$.

\item If, for some $c\in(0,\infty)$, $r_N\to c^{-2}$, then
\begin{equation}\label{eq:crit-crossover-limit}
 \frac{M_N}{N^{3/4}}\Longrightarrow Y_c,
\end{equation}
where $Y_c$ has density \eqref{eq:pc-def}.
Equivalently, if $p_NN^{3/4}\to\rho\in(0,\infty)$, then $c=\rho^{-2}$ and the quadratic coefficient is $1/(48\rho^2)$.
\end{enumerate}
\end{corollary}

\begin{proof}
In the first regime Proposition~\ref{prop:annealed-critical-regimes}(1) gives $\theta_N\to0$.  Apply Theorem~\ref{thm:crossover-K} with $c=0$ to obtain \eqref{eq:crit-quartic-limit}; combining the same theorem with the explicit bound following from \eqref{eq:theta-full-expansion} gives \eqref{eq:strong-dense-critical-rate}.  In the second regime Proposition~\ref{prop:annealed-critical-regimes}(2) gives $\theta_N\to c/24$, and Theorem~\ref{thm:crossover-K} yields \eqref{eq:crit-crossover-limit}.  No separate Laplace argument is required.
\end{proof}

\appendix
\section{Optimal derivative bound and the double-well transition}\label{app:optimal-K1}

The main theorem uses only Proposition~\ref{prop:crossover-simple-factors}.  The purpose of this appendix is different: it determines the optimal uniform bound on the first derivatives of indicator Stein solutions for the crossover family and identifies the point at which the double-well geometry first worsens that optimal bound.  None of the arguments in the main text depends on this sharpening.

\subsection{Exact variational representation}

For fixed $x$, differentiation of \eqref{eq:fz-crossover} gives
\[
\sup_z|f_{c,z}'(x)|=\max\{Q_c(x)|\ell_c'(x)|,F_c(x)|r_c'(x)|\}.
\]
On $x\ge0$ one has $r_c'(x)<0$.  Indeed, this is immediate while $g_c(x)\le0$; once $g_c>0$ and $g_c'>0$, the standard Mills inequality $r_c<1/g_c$ gives $g_cr_c-1<0$.  Put $A_c(x)=Q_c(x)\ell_c'(x)$ and $B_c(x)=F_c(x)(-r_c'(x))$.
The identities \eqref{eq:ell-r-ode} give $A_c(x)+B_c(x)=1$.
Since $B_c>0$, $A_c<1$.  If $A_c\ge0$, both contributions to the pointwise derivative bound are at most $1$; if $A_c<0$, then $B_c=1-A_c>1$ and dominates.  As the tail contribution tends to $1$, the optimal uniform derivative bound is therefore
\[
 C_{K,1}(c)
 =1+\left[-\inf_{x\ge0}Q_c(x)\ell_c'(x)\right]_+.
\]
In particular, $C_{K,1}(c)=1$ if and only if $\ell_c'(x)\ge0$ for all $x\ge0$.

\subsection{Reduction to one distinguished point}

For the threshold analysis assume now $c>0$.  The case $c=0$ is already covered by the quartic argument above.  The sign question reduces to a single point.  Put $u_c(x)=\ell_c'(x)=1+g_c(x)\ell_c(x)$, $b_c=\sqrt{c/24}$ and $a_c=\sqrt{c/8}$.
On $(0,a_c)$ the drift is negative and $u_c'(x)=g_c'(x)\ell_c(x)+g_c(x)u_c(x)$.
At a zero $u_c(x_0)=0$ we have $\ell_c(x_0)=-1/g_c(x_0)$ and thus $u_c'(x_0)=-g_c'(x_0)/g_c(x_0)$.
Thus any zero to the left of $b_c$ is crossed downwards, whereas any zero to the right of $b_c$ is crossed upwards.  At a zero exactly at $b_c$, one has $u_c'(b_c)=0$ and $u_c''(b_c)=g_c''(b_c)\ell_c(b_c)=2b_c\ell_c(b_c)>0$, so this zero is a strict local minimum.  Since $u_c(0)=u_c(a_c)=1$ and $u_c>1$ for $x>a_c$, the zero-crossing directions (including the tangential case at $b_c$) show that non-negativity of $u_c$ on $[0,\infty)$ is equivalent to $u_c(b_c)\ge0$.
At $b=b_c$ one has $g_c(b)=-(2/3)b^3$, so the threshold equation is
\[
 \frac{F_c(\sqrt{c/24})}{p_c(\sqrt{c/24})}
 =\frac{3}{2(c/24)^{3/2}}.
\]

\subsection{Scalar threshold equation and uniqueness}

\begin{figure}[ht]
 \centering
 \includegraphics[width=0.78\textwidth]{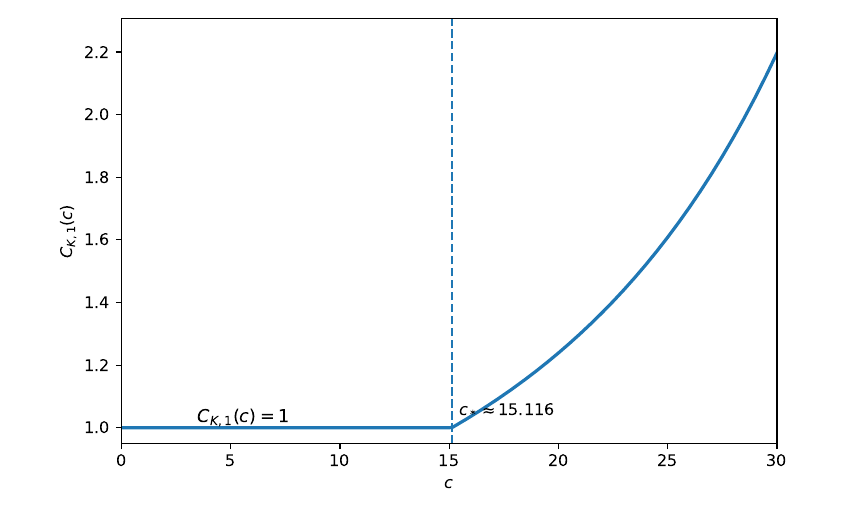}
 \caption{Numerical evaluation of the exact variational formula for the optimal uniform bound on the first derivatives of indicator Stein solutions $C_{K,1}(c)$. The optimal bound remains equal to $1$ for $0\le c\le c_*$, where $c_*\approx15.116$, and becomes strictly larger than $1$ beyond this threshold. The numerical evaluation and preparation of this figure were assisted by ChatGPT (OpenAI); the underlying formula and numerical output were independently verified by the author.}
 \label{fig:CK1-transition}
\end{figure}

This threshold equation has a unique solution.  To see this analytically, write $c=24b^2$, $y=b^4=c^2/576$ and $\psi(u)=-u^4/12+u^3/3-2u/3$.
The change of variables $t=b(1-u)$ gives $\ell_{24b^2}(b)=b\int_0^\infty e^{y\psi(u)}\,du$.
Hence the threshold condition is $H(y)=1$, where
\begin{equation}\label{eq:H-y}
 H(y)=\frac23y\int_0^\infty e^{y\psi(u)}\,du.
\end{equation}
Since $\psi'(u)+2/3=u^2(1-u/3)$ and $y\int_0^\infty\psi'(u)e^{y\psi(u)}\,du=-1$, integration by parts gives
\[
 H(y)-1=y\int_0^\infty u^2\left(1-\frac u3\right)e^{y\psi(u)}\,du.
\]
At a zero $H(y)=1$, normalize $e^{y\psi(u)}du$ to a probability measure and denote the corresponding random variable by $U$.  The preceding identity gives $\E U^3=3\E U^2$.  Integrating the derivative of $u e^{y\psi(u)}$ and using this relation yields
\[
 \frac{H'(y)}{H(y)}
 =\frac14\left\{\E(U-1)^2-1+\frac3y\right\}
 \qquad\text{whenever }H(y)=1.
\]
Thus every zero in $0<y<3$ is crossed strictly upwards.  On the other hand, $\psi(u)\ge0$ for $2\le u\le3$, so $H(y)\ge2y/3$; hence there is no zero for $y\ge3$.  Since $\psi(u)\le C-u^4/24$ for $u\ge0$ with a fixed finite $C$, one has $H(y)=O(y^{3/4})\to0$ as $y\downarrow0$.  Continuity, the strict upward crossing at each zero in $(0,3)$ and the absence of zeros for $y\ge3$ therefore give exactly one zero $y_*$.
Consequently there is a unique threshold $c_*=24\sqrt{y_*}$, with
$C_{K,1}(c)=1$ for $0\le c\le c_*$ and $C_{K,1}(c)>1$ for $c>c_*$.
Numerical quadrature of the one-dimensional integral \eqref{eq:H-y} gives $c_*\approx15.115657272651$.
The numerical value is included only for orientation; no certified numerical enclosure is used anywhere in the argument.
Figure~\ref{fig:CK1-transition} illustrates the resulting transition of the optimal uniform bound on the first derivatives of the indicator Stein solutions.

This threshold has no effect on the validity or on the $N^{-1/2}$ rate of the exchangeable-pair argument in the main text.  It only marks the point at which the \emph{optimal} uniform bound on the first derivatives of indicator Stein solutions ceases to equal the quartic value $1$.  This is why the main proof uses the simpler global estimate \eqref{eq:simple-K1}.

\FloatBarrier

\section*{Statements and Declarations}
\paragraph{Funding.}
No specific funding was received for this work.

\paragraph{Competing interests.}
The author has no relevant financial or non-financial interests to disclose.

\paragraph{Data availability.}
No datasets were generated or analysed during the current study.

\paragraph{Author contributions.}
Peter Eichelsbacher is the sole author and is responsible for all aspects of the work.

\paragraph{AI assistance.}
The author used ChatGPT (OpenAI) as an assistive tool for language editing, literature cross-checking,
and numerical assistance for the computation reported in Appendix~A.3, including the displayed
value of $c_*$ and the corresponding numerical figure.  All mathematical arguments, references,
numerical values and final formulations were independently checked by the author, who takes full
responsibility for the content of the manuscript.

\end{document}